\documentclass[11pt]{article}

\usepackage{fullwidth} 
\usepackage{color}
\usepackage[color,matrix,arrow]{xy}
\usepackage[top=1.5in, bottom=1.5in, left=1in, right=1in]{geometry}
\usepackage{amsgen}
\usepackage{amsmath}
\usepackage{amstext}
\usepackage{amsbsy}
\usepackage{amsopn}
\usepackage{amsfonts}
\usepackage{amssymb}
\usepackage{eepic}
\usepackage{graphicx}
\usepackage{tikz}
\usepackage{epsf}
\usepackage{pstricks}
\usepackage{mathrsfs}
\usepackage{faktor}
\usepackage{amsthm}
\xyoption{all}

\def\Box{\square}
\def\edge{\relbar\joinrel\relbar}
\def\longedge{\relbar\joinrel\relbar\joinrel\relbar\joinrel\relbar}
\def\mapright#1{\smash{\mathop{\longrightarrow}\limits^{#1}}}
\def\arr{\mapright{}}
\def\tra#1{\smash{\mathop{\mid\kern
-1pt\joinrel\relbar\joinrel\relbar}\limits^{*}_{#1}}}
\def\longtra#1{\smash{\mathop{\mid\kern
-1pt\joinrel\relbar\joinrel\relbar\joinrel\relbar}\limits^{*}_{#1}}}
\def\vlongtra#1{\smash{\mathop{\mid\kern
-1pt\joinrel\relbar\joinrel\relbar\joinrel\relbar\joinrel\relbar}\limits^{*}_{#1}}}
\def\vvlongtra#1{\smash{\mathop{\mid\kern
-1pt\joinrel\relbar\joinrel\relbar\joinrel\relbar\joinrel\relbar\joinrel\relbar}\limits^{*}_{#1}}}
\def\vvvlongtra#1{\smash{\mathop{\mid\kern
-1pt\joinrel\relbar\joinrel\relbar\joinrel\relbar\joinrel\relbar\joinrel\relbar\joinrel\relbar}\limits^{*}_{#1}}}
\def\etra#1{\smash{\mathop{\mid\kern
-1pt\joinrel\relbar\joinrel\relbar}\limits_{#1}}}
\def\mapleft#1{\smash{\mathop{\longleftarrow}\limits^{#1}}}
\def\longbar{\relbar\joinrel\relbar}
\def\vlongrightarrow{\relbar\joinrel\longrightarrow}
\def\vlongleftarrow{\longleftarrow\joinrel\relbar}
\def\vvlongrightarrow{\relbar\joinrel\vlongrightarrow}
\def\vvvlongrightarrow{\relbar\joinrel\vvlongrightarrow}
\def\vvvvlongrightarrow{\relbar\joinrel\vvvlongrightarrow}
\def\vvvvvlongrightarrow{\relbar\joinrel\vvvvlongrightarrow}
\def\vvvvvvlongrightarrow{\relbar\joinrel\vvvvvlongrightarrow}
\def\vvvvvvvlongrightarrow{\relbar\joinrel\vvvvvvlongrightarrow}
\def\vdashast{\stackrel{\ ,*}{\smash\vdash\vphantom{a}}}
\def\A{{\cal{A}}}
\def\X{{\cal{X}}}
\def\iff{\Leftrightarrow}
\def\Rw{\Rightarrow}
\def\oo{\overline}
\def\wt{\widetilde}
\def\wh{\widehat}
\def\B{{\cal{B}}}
\def\C{{\cal{C}}}
\def\F{{\cal{F}}}
\def\I{{\cal{I}}}
\def\L{{\cal{L}}}
\def\M{{\cal{M}}}
\def\N{\mathbb{N}}
\def\U{{\cal{U}}}
\def\S{{\cal{S}}}
\def\SB{\mathbb{SB}}
\def\pav{\mbox{Pav}}
\def\rat{\mbox{Rat}\ ,}
\def\acc{\mbox{acc}}
\def\aut{\mbox{Aut}}
\def\at{\mbox{At}}
\def\fct{\mbox{fct}}
\def\nbh{\mbox{nbh}}
\def\sji{\mbox{sji}}
\def\rec{\mbox{Rec}}
\def\dep{\mbox{dep}}
\def\dim{\mbox{dim}}
\def\fix{\mbox{Fix}}
\def\per{\mbox{Per}}
\def\irr{\mbox{Irr}\ ,}
\def\dom{\mbox{dom}}
\def\gcd{\mbox{gcd}}
\def\sol{\mbox{Sol}}
\def\lcp{\mbox{lcp}}
\def\lcm{\mbox{lcm}}
\def\reg{\mbox{Reg}\ ,}
\def\sing{\mbox{Sing}\ ,}
\def\pr{\mbox{Pr}}
\def\ker{\mbox{Ker}\ ,}
\def\pref{\mbox{Pref}\ ,}
\def\dw{\!\downarrow\ ,}
\def\upset{\!\uparrow\ ,}
\def\tr{\mbox{tr}}
\def\conj{\mbox{Conj}}
\def\equa{\mbox{Eq}}
\def\endo{\mbox{End}}
\def\gp{\mbox{Gp}}
\def\mt{\mbox{MT}}
\def\max{\mbox{max}}
\def\min{\mbox{min}}
\def\mindeg{\mbox{mindeg}}
\def\lat{\mbox{Lat}}
\def\LR{\mbox{LR}}
\def\cl{\mbox{Cl}}
\def\rk{\mbox{rk}}
\def\krk{\mbox{Krk}}
\def\prk{\mbox{prk}}
\def\geo{\mbox{Geo}}
\def\het{\mbox{ht}}
\def\sup{\mbox{sup}}
\def\sub{\mbox{Sub}}
\def\suff{\mbox{Suff}}
\def\ad{\mbox{Ad}}
\def\flatx{\mbox{Fl}}
\def\flats{\mbox{Fl\ ,}}
\def\P{{\cal{P}}}
\def\TR{{\cal{TR}}}
\def\R{{\cal{R}}}
\def\G{{\cal{G}}}
\def\H{{\mathcal{H}}}
\def\BoR{{\cal{BR}}}
\def\D{{\cal{D}}}
\def\V{{\cal{V}}}
\def\Zcal{{\cal{Z}}}
\def\Z{\mathbb{Z}}
\def\x{\!\cdot\!}
\def\p{\varphi}
\def\inv{^{-1}}
\def\Pnv{{\p \inv}}
\def\la{\langle}
\def\ra{\rangle}
\def\bow{\rtimes}
\def\bi{\begin{itemize}}
\def\ei{\end{itemize}}
\def\beq{\begin{equation}}
\def\eeq{\end{equation}}
\def\lexp#1{\ensuremath{{}^{#1}}}
\def\J{{\cal{J}}}
\def\AB#1{B(1,{#1})}
\def\lc{\left\lceil}   
\def\rc{\right\rceil}
\def\lf{\left\lfloor}   
\def\rf{\right\rfloor}
\def\SPC{\mbox{SPC}}
\def\un{\underline}
\def\GM{\operatorname{GM}}
\def\GGM{\operatorname{GGM}}
\def\LLM{\operatorname{LLM}}
\def\RLM{\operatorname{RLM}}
\def\RM{\operatorname{RM}}
\def\CM{\operatorname{CM}}
\def\LM{\operatorname{LM}}

\theoremstyle{plain}
\newtheorem{T}{Theorem}[section]
\newcommand{\bt}{\begin{T}}
\newcommand{\et}{\end{T}}
\newcommand{\ftd}{$\square$\end{T}}

\newtheorem{Proposition}[T]{Proposition}
\newcommand{\bp}{\begin{Proposition}}
\newcommand{\ep}{\end{Proposition}}
\newcommand{\fpd}{$\square$\end{Proposition}}

\newtheorem{Lemma}[T]{Lemma}
\newcommand{\bl}{\begin{Lemma}}
\newcommand{\el}{\end{Lemma}}
\newcommand{\fld}{$\square$\end{Lemma}}

\newtheorem{Corol}[T]{Corollary}
\newcommand{\bc}{\begin{Corol}}
\newcommand{\ec}{\end{Corol}}
\newcommand{\fcd}{$\square$\end{Corol}}

\newtheorem{Result}[T]{Result}
\newcommand{\br}{\begin{Result}}
\newcommand{\er}{\end{Result}}
\newcommand{\frd}{$\square$\end{Result}}

\theoremstyle{definition}
\newtheorem{Example}[T]{Example}
\newcommand{\be}{\begin{Example}}
\newcommand{\ee}{\end{Example}}

\newtheorem{Problem}[T]{Problem}
\newcommand{\bq}{\begin{Problem}}
\newcommand{\eq}{\end{Problem}}

\newtheorem{Remark}[T]{Remark}
\newcommand{\brm}{\begin{Remark}}
\newcommand{\erm}{\end{Remark}}

\newtheorem{Definition}[T]{Definition}
\newcommand{\bd}{\begin{Definition}}
\newcommand{\ed}{\end{Definition}}

\newtheorem{Construction}[T]{Construction}
\newcommand{\bco}{\begin{Construction}}
\newcommand{\eco}{\end{Construction}}

\newtheorem{Computation}[T]{Computation}
\newcommand{\bcomp}{\begin{Computation}}
\newcommand{\ecomp}{\end{Computation}}

\newcommand{\Qed}{\hfill
$\Box$
\par\bigbreak}

\normalbaselines
\def\abstract#1{\par\bigskip
\begingroup\small
\baselineskip=12truept
\begin{center}ABSTRACT\end{center}
\par\medskip\par\noindent
\null\hfill\hbox{\vbox{\hsize=5truein\noindent#1}}
\hfill\null\par\endgroup\par}
\def\up{^{{\rm up}}}
\def\tbpav{\mbox{TBPav}}
\def\bpav{\mbox{BPav}}

\title{Master List of Examples in Complexity Theory of Finite Semigroup Theory }

\author{Stuart Margolis and John Rhodes\footnote{We thank Jean-Camille Birget 
for writing the first draft of this document.}}
 
\date{\today}                 

\begin{document}
\maketitle

\tableofcontents

\pagebreak

\begin{center}

Abstract

\end{center}

This document gives a list of finite semigroups that are
interesting from the point of view of Krohn-Rhodes complexity theory. 
The list will be expanded and updated as ``time goes by''. \footnote{MSC classes: Primary 20M35, 20M07, 20M10, 20M20. \ Secondary 54H15.}



\section{Some background}

The purpose of this document is to maintain a list of examples that illustrate many aspects of the complexity theory (we write ``complexity'' for Krohn-Rhodes complexity) of finite semigroups. It is best for the reader to dive into the examples and learn the background material on the fly. Standard references for complexity theory are \cite[Chapters 7-9]{Arbib} and \cite{qtheory}. Important articles that give further background are \cite{Trans, FlowsI, complexity1, complexityn}. We have included a number of Appendices that covers much of the background material.

The document is not meant for publication. It will be updated periodically and maintained on the Research Gate sites of the authors and in the ArXiv. If you have comments, corrections, questions and/or examples you'd like to be included, please contact Stuart at: margolis@math.biu.ac.il. 

We'd also like to maintain a mailing list to let you know when the document is updated. If you are interested please send an email to Stuart at the above address.

\section{The Standard Format for Examples}

The following notation is standard and can mostly be immediately understood
by those doing research in finite semigroup theory.  
For all undefined notation see \cite{qtheory}, \cite{Arbib}, \cite{Eilenberg}. For further background see the Appendices to this paper and the references listed there.

The format will include the following.

\bigskip

\begin{enumerate}

\item{A regular Rees matrix semigroup $\mathcal{M}^{0}(G,A,B,C)$ with structure 
matrix $\ C: B \times A \to G^0 \ $ 
with $G$ a finite group, called the structure group  $\ G^0$ is $G$ with a zero added; $A$ and $B$ are finite 
non-empty sets. The matrix $C$ is {\em regular}, which means that every
row and every column has at least one non-0 entry. By Rees Theorem, finite regular Rees matrix semigroups are precisely the class of finite 0-simple semigroups.

Moreover, there are no proportional rows and columns in $C$. 
That is, $(\exists g \in G)(\forall a \in A)[g\ C(b_1,a) = C(b_2,a)]$
$\ \Rightarrow\ $ $b_1 = b_2\ $ (i.e., proportional rows are equal);
and dually on the right for columns. All 0-simple semigroups will have this property.
When $G$ is a non-trivial group this makes $\mathcal{M}^{0}(G,A,B,C)$  a group mapping semigroup denoted throughout as a $\GM$ semigroup. This means 
that the self-action of $\mathcal{M}^{0}(G,A,B,C)$ on the {\em left} and on the {\em right} of itself are faithful. For a 0-simple semigroup this is equivalent to saying that there are no proportional rows and columns in $C$. See Appendix \ref{Semiloc} for further details on $\GM$ semigroups and related material. 
} 

\item{Our semigroups $S$ will be subsemigroups of the translational hull 
$\Omega(I(S))$ of $I(S)=\mathcal{M}^{0}(G,A,B,C)$ that contain $I(S)$. 
$I(S)$ is then the unique 0-minimal ideal of $S$ and $S$ acts faithfully on the left and on the right of $I(S)$

More specifically, this is how we construct such semigroups. Let ${\rm RM}(B,G)$ be the semigroup of all {\em row-monomial} $B \times B$ 
matrices with coefficients in $G^0$; i.e., each row has at most one non-zero 
entry. Dually, ${\rm CM}(A,G)$ is the semigroup of all column-monomial 
$A \times A$ matrices with coefficients in $G^0$, where each column has at most
one non-zero entry. We note that ${\rm RM}(B,G)$ (${\rm CM}(A,G)$) is isomorphic to the (reverse) wreath product $G \wr PT_{R}(B)$ ($G \wr^{*} PT_{L}(A)$) of $G$ with the semigroup of all partial functions acting on the right of $B$ (left of $A$).

A matrix $X \in {\rm RM}(B,G)$ is said to be {\em linked} with a matrix 
$Y \in {\rm CM}(A,G^0)$ \ iff \ $CY = XC$. The collection of all such pairs is the translational hull $\Omega(I(S))$ of $I(S)$. Since we assume that $\mathcal{M}^{0}(G,A,B,C)$ is a $\GM$ 0-simple semigroup it follows that every element of $\rm {RM}(B,G)$ is linked with at most 1 element of ${\rm CM}(A,G)$ and vice versa.

Using the above notation, a $\GM$ semigroup gives a faithful transformation semigroup $(G \times B,S)$, where we identify a fixed 
$\mathcal{R}$-class in $I(S)$ with the set $G \times B$. $S$ induces an action on $B$ whose faithful image is called its {\em right letter mapping} image and denoted by $\RLM(S)$. We thus have the transformation semigroup $(B,\RLM(S))$. $(G \times B,S)$ can be proved to be a subsemigroup of $G \wr (B,\RLM(S))$. It is known that the computability of complexity can be reduced to deciding if given a $\GM$ semigroup $S$ whether its complexity $Sc$ satisfies $Sc = \RLM(S)c$ or $Sc = 1+ \RLM(S)c$. This is why $\GM$ semigroups play a central role in complexity theory and why all our examples will be $\GM$ semigroups. 
}

\item{
See Appendix \ref{Semiloc} for the basics of semilocal theory and the references \cite[Chapter 7-8]{Arbib} and \cite[Section 4.7]{qtheory} for more details on all of these assertions.}

\end{enumerate}


\section{The Master List}

\medskip
\noindent
Examples are given in the form of a $\GM$ Rees matrix semigroups and additional generators. For historic reasons, instead of $C$ as above we use the transpose
$C^{T}: A \times B \to G^0$. Generators are given (usually not a minimal set) of the form 
$(a,g,b) \in \mathcal{M}^{0}(G,A,B,C)$ that generate all of $\mathcal{M}^{0}(G,A,B,C)$, and a
collection of elements of ${\RM}(B,G)$, linked with some necessarily 
unique elements of ${\CM}(A,G)$. We write an element $X$ of ${\RM}(B,G)$ as a $G$-labeled partial function. The edge $i \rightarrow gj, i,j \in B, g \in G$ means that the $(i,j)$ entry of $X$ is $g$. Thus we identify $X$ with the set of such edges, one for each non-zero row of $X$. We have the dual representation for ${CM}(A,G)$. Given a generator $f \in \RM(B,G)$ we do not give the unique $f^{*} \in \CM(A,G)$ such that $(f,f^{*})$ is in the translational hull of $\mathcal{M}^{0}(G,A,B,C)$. The reader can compute that directly. 


\medskip
\noindent
With this notation, given a $\GM$ transformation semigroup $(G \times B,S)$ and $s \in S$ given as above by a $G$-labeled partial function, 
we have $(g,b)s=(gh,b')$ if $b \rightarrow hb'$ is an edge of $s$ and undefined if there is no edge beginning with $b$ in $s$. Thus $S$ acts faithfully as partial maps on $G \times B$. The image of the action induced by $(G \times B,S)$ on $B$ defines the $\RLM$-transformation semigroup, $(B,\RLM(S))$. As noted above, our purpose is to determine if $Sc = \RLM(S)c$ or $Sc = 1+ \RLM(S)c$. This is done by using the Theory of Flows. See Appendix \ref{Flows} for the definition and details. For further details, see the paper \cite{FlowsI}.

\medskip
\noindent
All  computations take place in the Evaluation Transformation Semigroup $\mathcal{E}(L)=(\operatorname{States},\operatorname{Eval}(S))$, where $L$ is the Rhodes lattice $Rh_{B}(G)$. See Appendix \ref{Eval} for the definitions and \cite{Trans} for more details. We do our computations for flows using the Rhodes lattice $Rh_{B}(G)$ rather than the Set-Partition lattice $SP(G \times B)$, as was done in Chapter 4 of \cite{qtheory} and Section 5 of \cite{Trans}. See Appendix \ref{SetPart} for the connection between these lattices. See Section \ref{SPCNotation} for the notation we use for elements of Rhodes lattices.
 
\smallskip

An Example includes the following information:

\begin{enumerate}

\item{The name of the example.}

\item{Description of the 0-minimal ideal ${\cal M}^0(G; A, B, C)$.}

\item{The set of generators.}

\item{A computation in the Evaluation Transformation Semigroup that constructs a flow of the appropriate complexity or proves that no such flow exists.See Appendix \ref{Flows} and Appendix \ref{Eval}. }

\item{Discussion of all of the above.}

\item{History of the example and references to the literature.}

\end{enumerate}

\bigskip

\begin{enumerate}

\item{{\bf The Tall Fork}\label{TF}

\medskip

The Tall Fork $TF$ is defined as follows:

\noindent Structure group: \ ${\mathbb Z}_2 = \{\pm 1\}$, so
$\ ,{\mathbb Z}_2^{ \ 0} = \{0,\pm 1\}$, written multiplicatively

\bigskip

\begin{tabular}{r||r|r|r|}
$\cdot$  & 0 & 1    & $-1$  \\ \hline \hline
0        & 0 & 0    & 0  \\ \hline
1        & 0 & 1    & $-1$ \\ \hline
-1       & 0 & $-1$ & 1  \\ \hline
\end{tabular}

\bigskip

\noindent The structure matrix of $I(TF)$ is given by 

\bigskip

$C^{T} \ = \ \ $
\begin{tabular}{r|c c|r r r r|}
      &$1'$&$3'$& 1 & 2 & 3 & 4  \\ \hline
$a_1$ & 1 & 1   & 0 & 0 & 0 & 0  \\ 
$a_2$ & 1 & 0   & 0 & 0 & 0 & 0  \\ 
$a_3$ & 0 & 1   & 0 & 0 & 0 & 0  \\ \hline 
$a_4$ & 0 & 0   & 1 & 1 & 0 & 0  \\ 
$a_5$ & 0 & 0   & 0 & 1 & 1 & 0  \\ 
$a_6$ & 0 & 0 & 0 & 0 & 1 & 1  \\ 
$a_7$ & 0 & 0  & 1 & 0 & 0 & 1  \\ \hline
\end{tabular}

\bigskip

\noindent The generators of $TF$ are all the elements of $I(TF)$ together with the following elements of ${\rm RM}(B,\mathbb{Z}_{2})$. 

$\sigma = (1' \ 3')$ \ \ (transposition of columns), so 
$\langle \sigma \rangle \simeq {\mathbb Z}_2$;

$\tau = (1 \ 2 \ 3 \ 4)$ (cyclic permutation of columns), so
$\langle \tau \rangle \simeq {\mathbb Z}_4$;

$r = \{(1' \to 1), \ (3' \to -3)\}$ \ \ \ (note the sign-change in 
column 3).

\pagebreak

\noindent The poset of $\mathcal{J}$-classes is as follows. This explains the name the name ``Tall Fork'':

\bigskip

\hspace{1,8cm} ${\mathbb Z}_2$ \hspace{1,5cm} ${\mathbb Z}_4$

\smallskip

\hspace{2,2cm} \verb|\       /| 

\smallskip

\hspace{2,4cm} $[{\mathbb Z}_2\ , r\ , {\mathbb Z}_4]$
 \hspace{2,0cm} (null class)

\smallskip

\hspace{3,1cm} $|$

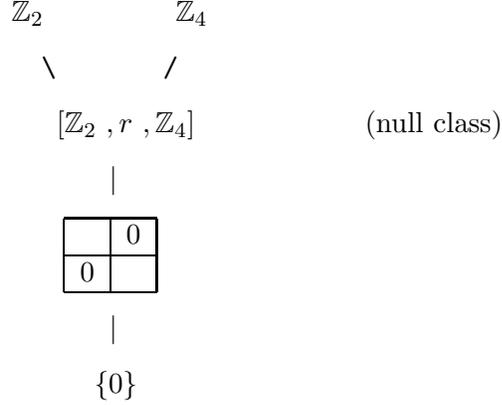
\smallskip

\hspace{2,5cm}
\begin{tabular}{|c|c|} \hline
  & 0 \\ \hline
0 &   \\ \hline
\end{tabular}

\smallskip

\hspace{3,1cm} $|$

\smallskip

\hspace{2,9cm} $\{0\}$

\bigskip

\bigskip

where the matrix above $\{0\}$ represents the non-zero $\mathcal{J}$-class of $\mathcal{M}^{0}(Z_{2},A,B,C)$.

\noindent Note that the longest chain of non-aperiodic $\mathcal{J}$-classes in $\RLM(TF)$ has length 1. It follows that $\RLM(TF)c = 1$ by the Depth Decomposition Theorem \cite{TilsonXII}. We show that $(TF)c = 2$, by showing that no aperiodic flow exists for $TF$. Proofs appear in \cite{KernelSystems}, Chapter 4 of 
\cite{qtheory} using the Presentation Lemma and \cite[Section 6.]{Trans} using flows over the Set-Flow lattice $SP(G \times B)$. We use the Rhodes lattice $Rh_{B}(G)$, but the justification is the same as the reference in \cite[Section 6.]{Trans}. We give the proof here for the convenience of the reader. See Appendix \ref{SetPart} for information on the Set-Partition Lattice $SP(G \times B)$ and the Rhodes Lattice
$Rh_{G}(B)$. See Appendix \ref{Eval} for the Evaluation Transformation Semigroup and Section \ref{SPCNotation} for the notation we use for elements of the Rhodes Lattice. Here is the computation.

\medskip

$(\{1'\})/\langle 1 \rangle )(\sigma)^{\omega + *} \ = \ \{1'  \ 3'\}/\langle 1 \ 1\rangle$. We used the Tie-Your-Shoe Lemma \cite[Lemma 6.1]{Trans} to justify that the result of applying $(\sigma)^{\omega + *}$ has one partition class.

\medskip

\medskip

$(\{1' \ 3'\}) /\langle 1 \ 1\rangle) r \ = \{1 \ 3\}/\langle 1 \ -1\rangle$

\medskip

$(\{1| \ 3\}/\langle 1 \ -1\rangle) (\tau)^{\omega + *}$ is not a cross-section and thus is the contradiction $\Rightarrow \Leftarrow$, the contradiction in the Rhodes lattice.

\medskip

\noindent We have used the Vacuum and the Tie-Your-Shoes Lemma to justify the claim that this last computation is correct. It follows that no aperiodic flow exists and thus $(TF)c = 2$.

It is not difficult to show that the maximal complexity of any submonoid of $TF$ is 1 \cite[Chapter 4]{qtheory}. 
Therefore the complexity function $c$ is not local, that is, it is not the maximum of the complexities over submonoids of a semigroup. In its role as a counterexample to locality, we have one of the reasons that $TF$ was defined.

\bigskip

\noindent History: From J.\ Rhodes (1977) , \cite{qtheory} and \cite{KernelSystems}.}

\bigskip

\item{{\bf  A subsemigroup of the Tall Fork}\label{TFA1} 

\smallskip
 
\noindent The definition of the semigroup $TFA1$ is the same as Example \ref{TF} EXCEPT that row $a_{1}$ is removed from the transpose of the structure matrix $C^{T}$ of the 0-minimal $\mathcal{J}$-class. This new example has $c = 1$, and here is an aperiodic 2-state flow. We first define an aperiodic transformation semigroup that $A = (\{p,q\},RZ(\{p,q\})^{1})$, where $RZ(\{p,q\}^{1}$ is the flip-flop, consisting of the two constant functions (which are right-zeroes and explain the letters $RZ$) and the identity. 

%


We define a function $F:\{p,q\} \rightarrow Rh_{B}(G)$ as follows:

\noindent

$pF = \{1' \ 3'\}/\langle1 \ -1 \rangle$ 
\noindent

$qF=\{1 \ 2 \ 3 \ 4\}/\langle 1 \ 1 \ 1 \ 1 \rangle$

%
%


\smallskip

\noindent Now we cover generators $x \in X$ of $TFA1$.

\begin{enumerate}

\item{If  $x = (a,g,b), b \in \{1, \ 2, \ 3, \ 4\}$, then we cover $x$ by is $C(q)$, the constant function to $q$.}

\item{If $x = (a,g,b), b \in \{1', \ 3'\}$ then we cover $x$ by $C(p)$, the constant function to $p$.}

\item{$\sigma$ and $\tau$ are covered by the identity function.}

\item{$r$ is covered by $C(q)$.}

\end{enumerate}

\noindent This defines an $X$-automaton structure on the states of the transformation semigroup $A$. We now verify that $F$ is a flow.

\begin{enumerate}
\item{If $x = (a, g, b), b \in \{1', \ 3'\}$, then $$(pF)x \leq \{b\}/\langle 1\rangle  \leq (px)F = pF =  \{1' \ 3'\}/\langle 1 \ -1\rangle $$ and 
$$(qF)x \leq {b}/\langle 1 \rangle  \leq (qx)F = pF = \{1', \ 3'\}/\langle 1 \ -1 \rangle $$.}

\item{If $x = (a, g, b), b \in \{1, \ 2, \ 3, \ 4\}$, then $$(pF)x \leq \{b\}/\langle 1 \rangle \leq (px)F = qF = \{1 \ 2 \ 3 \ 4\}/\langle 1 \ 1 \ 1 \ 1 \rangle $$ and 
$$(qF)x \leq {b}/\langle 1 \rangle  \leq (qx)F = qF = \{1 \ 2 \ 3 \ 4\}/\langle 1 \ 1 \ 1 \ 1\rangle $$.}

\item{If $x = \sigma$, then $$(pF)x = (\{1' \ 3'\}/\langle 1 \ -1 \rangle )\sigma = \{1 \ 3\}/\langle 1 \ 1 \rangle  \leq (px)F $$ and 
$$(qF)x  = (\{1 \ 2 \ 3 \ 4\}/ \langle 1 \ 1 \ 1 \ 1 \rangle )\sigma = \emptyset \leq (qx)F $$}

\item{If $x = \tau$, a computation similar to that for $\sigma$ shows that the conditions verifying that $F$ is a flow hold.}

\item{If $x = r$, then $$(qF)x = (\{1 \ 2 \ 3 \ 4\}/\langle 1 \ 1 \ 1 \ 1\rangle )x = \emptyset \leq (qx)F$$ and 
$$(pF)x = (\{1' \ 3'\}/\langle 1 \ -1\rangle )x = \{1 \  3\}/\langle 1 \ 1\rangle \leq \{1 \ 2 \ 3 \ 4\}/\langle 1 \ 1 \ 1 \ 1\rangle  = (px)F$$}

\end{enumerate}

This completes the proof that $F$ is a flow. We have done the computation in full detail in this first check that a given function defines a flow and hope that helped the reader. In the rest of the examples, we'll give shorter versions of these computations.

\noindent History: By John Rhodes 1977, \cite{KernelSystems}}

\bigskip

\item{{\bf Using the Vacuum}

We saw that removing row $a_1$ from the structure matrix in Example \ref{TF} dropped the complexity from 2 to 1. See Example \ref{TFA1}. In this Example we keep the structure matrix from Example \ref{TFA1} but add a new generator. We'll see that this causes the complexity to rise to 2 once again. This Example is designed to show the effect of using the Vacuum operator $V$. See Definition \ref{Flowops} paragraph \ref{Vac}.

\medskip

\noindent Structure group with zero: 
$ \ {\mathbb Z}_2^{ \ 0} = (\{0, \pm1\}, \cdot)$,

\noindent with structure matrix 

\bigskip

$C^{T} \ = \ $
\begin{tabular}{r|c c|r r r r|}
      &$1'$&$3'$& 1 & 2 & 3 & 4  \\ \hline
$a_2$ & 1 & 0   & 0 & 0 & 0 & 0  \\
$a_3$ & 0 & 1   & 0 & 0 & 0 & 0  \\ \hline
$a_4$ & 0 & 0   & 1 & 1 & 0 & 0  \\
$a_5$ & 0 & 0   & 0 & 1 & 1 & 0  \\
$a_6$ & 0 & 0 & 0 & 0 & 1 & 1  \\
$a_7$ & 0 & 0  & 1 & 0 & 0 & 1  \\ \hline
\end{tabular}
%
%

\bigskip

\noindent This matrix is the structure matrix from Example \ref{TFA1}. The semigroup $UTV$ (Using the Vacuum) has generators all the elements in $I(UTV)$ plus the following.

\medskip

$r: \ 1' \mapsto 1, \ 3' \mapsto -\!3$

\medskip

$t: \ 1' \mapsto 1, \ 3' \mapsto 3$

\medskip

$\sigma = (1' \ 3')$, \ \ \   \ \ \ $\tau = (1 \ 2 \ 3 \ 4)$,


We have added the element $t$ to the generators of Example \ref{TF} modulo using the structure matrix from Example \ref{TFA1}.

\medskip

\noindent We claim that $(UTV)c = 2$. Here are the computations.

\medskip
\noindent

$\{1'\}/\langle 1 \rangle (\sigma)^{\omega + *} =\{1'\ \mid \ 3'\}/\langle 1|1\rangle $ and $(\{1'\ | \ 3'\}/\langle 1\mid 1\rangle)t = \{1\ \mid \ 3\}/\langle 1\mid 1\rangle$.

\medskip
\noindent 

Now $\{1 \ \mid \ 3\}/\langle 1|1\rangle (\tau)^{\omega+*} = (\{1\  2 \ 3\  4\}/\langle 1 \ 1 \ 1 \ 1\rangle )$. 

\medskip
\noindent

Therefore $(\{1\ \mid \ 3\}/\langle 1\mid 1\rangle )V= \{1  3\}/\langle 1 \ 1\rangle $. That is, the Vacuum $V$ forces $1$ and $3$ to be in one partition class.

\medskip
\noindent

We now have that $(\{1'\ | \ 3'\}/\langle 1|1 \rangle )t = \{1\ | \ 3\}/\langle 1 \ 1 \rangle $ and this forces $(\{1'\ | \ 3'\}/\langle 1|1 \rangle )V = \{1'\   \ 3'\}/\langle 1 \ 1\rangle $. But then $(\{1' \ 3'\}/\langle 1 \ 1 \rangle )r = \{1 \ 3\}/\langle 1 \ -1 \rangle $. 

\medskip
\noindent

This last $SPC$ is precisely the 
State in Example \ref{TF} that lead to a contradiction. Therefore, $(UTV)c=2$.

\noindent History: John Rhodes, June 2019.} 

\bigskip

\item{{\bf BIRIP (Back Inject Rest in Peace)}\label{BIRIP}

(very important example from Karsten Henckell 2006)

\noindent We first define, for every $k \geq 1$, the $2k \times k$ matrix $M_k$ as follows.

\bigskip

$M_k \ = \ \ \ $
\begin{tabular}{r|c c c c c|}
      & 1 & 2 & 3 & \ldots & k  \\ \hline
$a_1$ & 1 & 1 & 0 & \ldots & 0  \\
$a_2$ & 0 & 1 & 1 & \ldots & 0  \\
 .    &   &   & cycled &   & \\ 
 .    &   &   & \ldots &   &\\
$a_k$ & 1 & 0 & \ldots & 0 & 1 \\ \hline
$a_1'$&  &  &  &  & \\ 
 .    &   &   &  & &   \\
 .    &   &   &{\bf I}$_k$ & & \\ 
 .    &   &   &  &   &\\
$a_k'$ &  &  &  &  & \\ \hline
\end{tabular}

\bigskip

\bigskip

\noindent Maximal subgroup with zero: \ ${\mathbb Z}_2^{ \ 0}$.

\bigskip

\noindent The structure matrix is 

\bigskip

\begin{tabular}{|c c|c c c c|} 
$1'$&$3'$& 1 & 2 &3 & 4  \\ \hline
 &  & & & &  \\
 & $M_2$ & & {\bf 0} & & \\
 &   &   & & &   \\ \hline
 &   &   & & &   \\
 &   &   & & &   \\
 &   &   & & &   \\
 & {\bf 0} &  & $M_4$  & &   \\
  &   &   &  &  &  \\
  &   &   & & & \\       
 &  &  &  & & \\ \hline
\end{tabular}

\bigskip

\noindent with generators I(BIRIP), plus

\smallskip

$g_2 = (1' \ 3')$ 

\smallskip

$a: \ 1 \mapsto 4, \ 2 \mapsto 3$

\smallskip

$s: \ 1' \mapsto 1, \ 3' \mapsto -\!3$ \ (negative).

\medskip

\noindent Then $(BIRIP)c = 1$. The reader is invited to prove that $RLM(BIRIP)$ has complexity 1 and that there is a 3-state flow with semigroup $RZ(3)^{1}$ proving that $(BIRIP)c = 1$. You can also look at Section 3.5 of \cite{complexity1} for details.

 \item{{\bf CBIRIP (compact back inject rest in peace)}\label{CBIRIP}

\bigskip

\begin{tabular}{|c c c c|}
1 & 2 & 3 & 4  \\ \hline
${\mathbb Z}_2^{ \ 0}$  & & &  \\ \hline 
 & & &   \\ 
 & & &   \\
 & $M_4$ & &   \\
 & & &  \\ \hline
\end{tabular}

\bigskip

\noindent with generators at I(CBIRIP), plus

\smallskip

$a: \ 1 \mapsto 4, \ 2 \mapsto 3$

$b = (1 \ 3)\ ,(2 \ 4)$

$s: \ 2 \mapsto 1, \ 2 \mapsto -\!3$

\medskip

\noindent Then $c = 1$ with a flow defined on $RZ(2)^1$. This is similar to Example \ref{BIRIP}. See Example 5.11 of \cite{deg2part2} for details.}
\bigskip

 \item{{\bf Geometry of \ ${\mathbb Z}_2 \times {\mathbb Z}_2$ 
$+$ ``Rube Goldberg''}

We define the semigroup $RG1$ with maximal subgroup ${\mathbb Z}_2$. The structure matrix for $I(RG1)$ is defined as follows.

\medskip

\bigskip

\begin{tabular}{r|c c|c c c c| c c|}
      &$1'$&$3'$& 1 & 2 & 3 & 4 &$1^R$&$2^R$  \\ \hline
$a_1$ & 1 & 1     & 0 & 0 & 0 & 0 & 0 & 0 \\
$a_2$ & 1 & 0     & 0 & 0 & 0 & 0 & 0 & 0 \\
$a_3$ & 0 & 1     & 0 & 0 & 0 & 0 & 0 & 0 \\ \hline
$a_4$ & 0 & 0     & 1 & 0 & 0 & 0 & 0 & 0 \\
$a_5$ & 0 & 0     & 0 & 1 & 0 & 0 & 0 & 0 \\
$a_6$ & 0 & 0     & 0 & 0 & 1 & 0 & 0 & 0 \\
$a_7$ & 0 & 0     & 0 & 0 & 0 & 1 & 0 & 0 \\ \hline
$a_8$ & 0 & 0     & 0 & 0 & 0 & 0 & 1 & 0 \\ 
$a_9$ & 0 & 0     & 0 & 0 & 0 & 0 & 0 & 1 \\ \hline
\end{tabular}

\bigskip

Omitting row $a_1$ gives the $8 \times 8$ identity matrix. Thus the only elements in $I(RG(2))$ that don't act by partial bijections on 
${\mathbb Z}_{2} \times B$ have $A$ coordinate equal to $a_{1}$. All the other generators will act by partial bijections on 
${\mathbb Z}_{2} \times B$ , so $RG1$ is ``almost'' an inverse semigroup.

\bigskip
 
\noindent Generators: All the elements in $I(RG1)$ plus the following elements.
\smallskip

$a = (1\ 4)(2\ 3)$

\smallskip

$b = (1\ 3)(2\ 4)$

\smallskip

Note that the subsemigroup generated by $a$ and $b$ is isomorphic to the group ${\mathbb Z}_{2} \times {\mathbb Z}_{2}$ and acts transitively on $\{1\ , 2\ , 3\ , 4\}$.

\smallskip

$c= (1'\ 3')$

$s: \ 1' \mapsto 1, \ 3' \mapsto -\!3$

\smallskip

$r_{1}: 1 \mapsto 1^R, \ 2 \mapsto 2^R$

\smallskip

$x_{R}:1^{R} \to -2^{R}, \ 2^{R} \to 1^{R}$

\bigskip

We note that the 0-minimal ideal of $\RLM(RG1)$ is aperiodic  and all elements outside the 0-minimal ideal act by partial bijections. Therefore the Rees quotient of $\RLM(RG1)$ by its 0-minimal ideal embeds in an appropriate symmetric inverse semigroup. It follows by the 
Fundamental Lemma of Complexity \cite{qtheory} that $\RLM(RG1)c=1$. 

We now define a flow from the transformation semigroup having the semigroup 
$RZ(3)^{1}$ acting on the following three States of the Evaluation Transformation Semigroup. That is, we have three states acted upon by $RZ(3)^{1}$ and we indicate the three elements of the Evaluation Transformation Semigroup to which they are mapped by the flow.

\begin{enumerate} 

\item{$\{1'\  3'\}/\langle 1 \  1 \rangle $}

\item{$\{1\ 3\ \mid \ 2 \ 4\}/\langle 1 \ -1 \ \mid \ 1 \ -1 \rangle  $}

\item{$\{1^{R}\ \mid \ 2^{R}\}/\langle 1\ \mid \ 1 \rangle $}

\end{enumerate}

We outline the computations proving that we have defined a flow. We cover $a,b,c$ and $x_{R}$ by the identity of $RZ(3)^{1}$. Since each of these elements is either not defined on one of the 3 states above or maps it to itself, we are done in this case. We cover $s$ by the constant map to state (b) and $r_1$ by the constant map to state (c). Finally we cover $(a,g,b) \in I(RG1)$ by the constant map to the unique state above whose set contains $b$. We leave it to the reader to show that this proves that we have defined a flow.
Therefore we have that $(RG1)c = 1$ and we have the following decomposition.
\smallskip

\begin{center}

 \ \ \ $RG1 \ \prec \ $
$({\mathbb Z}_2 \x {\mathbb Z}_2) \wr {RZ}(3)^1$
$\times$ ${\rm \RLM(RG1)}$ 

\end{center}

\medskip
\noindent

{\bf Activating Rube Goldberg's Contraption Changes the Complexity to 2}

We modify the generators of $RG1$ by keeping all the generators of $RG1$ except for $r_{1}$. We replace $r_{1}$ with $r_{2}$ that sends
$1$ to $1_{R}$ and $3$ to $-2_{R}$. Call the resulting semigroup $RG2$. The same computation for $RG1$ shows that $\RLM(RG)c = 1$.

Note now that $(\{1\ 3\ \mid \ 2 \ 4\}/\langle 1 \ -1 \ \mid \ 1 \ -1 \rangle )r_{2} = \{1_{R} \ 2_{R}\}/\langle 1 \ 1\rangle$. Acting on this
 state by $x_{R}$ sends it to $\{1_{R} \ 2_{R}\}/\langle 1 \ -1\rangle$. These two states are then joined by acting by $x_{R}^{\omega + *}$.
 But this join is the contradiction of the Rhodes lattice. Therefore $(RG2)c = 2$.

\bigskip
\noindent

$RG1$ and $RG2$ explain the use of the name ``Rube Goldberg''. Rube Goldberg was an American author who invented very complicated contraptions to do simple tasks. See here for a typical example of a Rube Goldberg contraption.
\medskip
\begin{center}
https://miro.medium.com/v2/resize:fit:828/format:webp/1*UXeSA9hDKfUEfw4cE2pIpA.png
\end{center}

In order to get the contraption to work, all the various pieces had to be activated. In our example, the element $x_{R}$ is a switch in our Rube Goldberg semigroup. In example $RG1$, the switch was not activated and the complexity of $RG1$ is 1. On the other hand, the element $r_{2}$ activates the switch by forcing the two elements $1_{R}$ and $2_{R}$ to be in one partition block. This leads the contraption to kick the complexity up to 2.

\noindent History: J.\ Rhodes, 11 August 2024}

\item{{\bf Example where an Aperiodic Flow Exists, but the Complexity is 2}\label{NoAFlow}

\medskip
\noindent
The examples up to now have computed an aperiodic flow in order to show that a $\GM$ semigroup $S$ such that $\RLM(S)c = 1$ satisfies $Sc = 1$. Here we give an example of a $\GM$ semigroup $S$ that admits an aperiodic flow, but $Sc=2$ because $\RLM(S)c=2$. We want to emphasize that the existence of an aperiodic flow only implies that $Sc = \RLM(S)c$ and thus there is a need to compute the complexity of $\RLM(S)$ as well. We remind the reader that if $S$ is a $\GM$ semigroup, then the cardinality of $\RLM(S)$ is strictly less than that of $S$, so inductive proofs are available.

We start with a general construction. Let $S$ be a $\GGM$ semigroup with $I(S) = \mathcal{M}^{0}(G,A,B,C)$. Since we are assuming that $S$ is 
$\GGM$ it is possible that $G$ is the trivial group. Let $H$ be a non-trivial group. Then the minimal ideal $K(H \times S)$ of 
$H \times S$ satisfies $K(H \times S) = H \times \{0\}$. We call the Rees quotient $H \times S/K(H\times S)$ the reduced direct product of $H$ and $S$ and denote it by $H \times_{r} S$. It is easy to show that $I(H \times_{r} S) = \mathcal{M}^{0}(H \times G,A,B,C)$ and that $H \times_{r} S$ is a $\GM$ semigroup. 

We now construct our example. Let $I_{4}= \mathcal{M}^{0}(1,\{1,\ldots, 8\}, \{1, \ldots, 4\}, M_{4}^{T})$ be the $\GGM$ 0-simple semigroup over the trivial group. Here $M_{4}$ is the matrix defined in Example \ref{BIRIP} in the case $k=4$.
By Lemma 4.1 of \cite{cremona}, the translational hull $\Omega(I_{4})$ of $I_{4}$ is isomorphic to the monoid of continuous partial functions on a $4$-cycle $\Gamma_{4}=(V,E)$ where $V =\{1, \ldots, 4\}$ and $E =\{\{i,i+1\} \mid i = 1 \ldots 4\}$ where we take the indices modulo $4$. This is the collection of all partial functions $f:\{1, \ldots, 4\} \rightarrow \{1, \ldots, 4\}$ such that if $x \in (V \cup E)$
is a vertex or edge of $\Gamma_{4}$, then either $xf^{-1}$ is empty or belongs to $V \cup E$.  This makes it fairly easy to check if a given partial function $f:\{1, \ldots, 4\} \rightarrow \{1, \ldots, 4\}$ belongs to the translational hull of $I_{4}$. See Example 5.4 of \cite{deg2part2} for the Green-Rees structure of $\Omega(I_{4})$.

Let $G$ be a non-trivial group and let $T_{4}(G) = G \times_{r} \Omega(I_{4})$ be the reduced direct product of $G$ with $\Omega(I_{4})$.
Then the 0-minimal ideal of $T_{4}$ is $\mathcal{M}^{0}(G,\{1,\ldots, 8\}, \{1, \ldots, 4\}, M_{4}^{T})$. Furthermore $T_{4}(G)$ is a $\GM$ semigroup and its $\RLM$ image is easily seen to be $\Omega(I_{4})$. 

It follows from the third item in Theorem \ref{typeIIprop} that the type II subsemigroup of $T_{4}(G)$ is isomorphic to the type II subsemigroup of $\Omega(I_{4})$. In particular the intersection of the type II subsemigroup of $T_{4}(G)$ with $I(T_{4}(G))$ is aperiodic. It now follows from Theorem 5.2 of \cite{FlowsI} that $T_{4}(G)$ admits a flow over the trivial transformation semigroup.  We have that the complexity of $T_{4}(G)$ is the same as its $\RLM$ image, $\Omega(I_{4})$. Example 5.4 of \cite{deg2part2} shows that $\Omega(I_{4})c=2$. Therefore $T_{4}(G)$ admits an aperiodic (even trivial) flow, but does not have complexity 1.

History: Stuart Margolis and John Rhodes, 2023 \cite{deg2part2}. }

\item{{\bf Example of a Semigroup with a Flow of Complexity 1 but no Aperiodic Flow}

Let $G$ be a group. Define $I_{4}(G) = \mathcal{M}^{0}(G,\{1,\ldots, 8\}, \{1, \ldots, 4\}, M_{4}^{T})$, as in Example \ref{NoAFlow}. We define $S_{4}(G)$ to be the translational hull of $I_{4}(G)$. We write $S_{4}$ for $S_{4}(1)$. As mentioned in Example \ref{NoAFlow}, $S_{4}$ is the monoid of continuous partial functions on a 4-cycle \cite{cremona} , \cite{deg2part2}. The reader can verify that $\RLM(S_{4}(G)) = S_{4}$. Example 5.4 \cite{deg2part2}, we have that $S_{4}c=2$. On the other hand, it is known that if $(g,b) \in G \times B, f \in S_{4}(G)$, then 
$|(g,b)f^{-1}\mid  \leq 2$. That is, every fiber of an element of $S_{4}$ has cardinality at most 2 in its action on $G \times B$. By definition this means that the transformation semigroup $(G \times B, S_{4}(G))$ has degree at most 2.  It follows from \cite{Tilsonnumber} , \cite{cremona} , \cite{deg2part2} that $S_{4}c \leq 2$.

Therefore, $S_{4}c = \RLM(S_{4})c = 2$. It follows from Theorem \ref{uniform} that $S_{4}$ admits a flow from a transformation semigroup of complexity of at most 1. We define a flow of complexity 1 and show that there is no aperiodic flow. The proof that a  transformation semigroup $(Q,S)$ such that every element in $S$ has degree at most $n$ satisfies $Sc \leq n$ \cite{Tilsonnumber} computes a flow of complexity at most $n-1$ in modern language. 

We recall the definition of the fiber graph of a transformation semigroup $X=(Q,S)$ of degree at most 2 from \cite{cremona} , \cite{deg2part2}. It is the graph $\Gamma(X)$ with vertices $Q$ and edges $\{u,v\}$ where $\{u,v\}$ is a fiber of cardinality 2 of some $s \in S$. Let $V(G) = (G \times \{1, \ 2, \ 3, \ 4\})$. Then for the transformation semigroup $(V(G), S_{4}(G))$, its fiber graph is a disjoint copy of $|G|$ copies of a 4-cycle. For each $g \in G$ there is a 4-cycle on the vertices $\{(g,1), \ (g,2), \ (g,3), \ (g,4) \}$. This can be verified by direct calculation or by Theorem 5.10 of \cite{cremona}. 

 Let $\mathcal{A}(G)$ be the set of anti-cliques (also known as independent sets) in the vertex set of the fiber graph of  $S_{4}(G)$. Theorem 3.5 of \cite{cremona} defines a reverse semidirect product 
$SIS(\mathcal{A}(G)) \rtimes 2^{\mathcal{A}(G)}$  of the symmetric inverse semigroup on $\mathcal{A}(G)$ with the semilattice of subsets of $\mathcal{A}(G)$. This is a semigroup of complexity 1. An action of this semigroup  on $\mathcal{A}(G)$ is defined. A flow to $(V(G), S_{4}(G))$ that sends each $X \in \mathcal{A}(G)$ to the $SPC$ whose partition is the singelton partition on $X$ with its unique cross section.  See Section 3 of \cite{cremona} for details.

We now show that $S_{4}(G)$ admits no aperiodic flow. In the next computation in the Evaluation ts of $S_{4}(G)$ the reader can verify directly that all elements belong to $S_{4}(G)$. Alternatively one can use the results of Section 5 of \cite{cremona} especially Proposition 5.6.

%
%
%
%

\medskip
\noindent

$(\{1\}/\langle 1 \rangle )(13)^{\omega+*} = \{1 \ \mid \ 3\}/\langle 1 \mid 1\rangle$. We then have 
$(\{1 \ \mid \ 3\}/\langle 1 \mid 1\rangle)(1 \  2 \ 3 \ 4)^{\omega + *} = \{1\  2 \ 3\  4\}/\langle 1 \ 2 \ 3 \ 4 \rangle $. By applying the vacuum $V$, we have that $(\{1 \ \mid  \ 3\}/\langle 1 \mid 1\rangle)V = \{1 \  3\}/\langle 1 \ 1\rangle$. 

\medskip
\noindent

Let $f = 1 \to 3 \ , 3 \to g1$ where $ g \neq 1 \in G$. Then $(\{1 \  3\}/\langle 1 \ 1\rangle)f = \{1 \  3\}/\langle 1 \ g\rangle$. Since 
$(1 \  2 \ 3 \ 4) \in S_{4}(G)$, we are led to a contradiction as in Example \ref{TF}. Therefore, $S_{4}(G)$ does not admit an aperiodic flow.

\medskip
\noindent
History: Stuart Margolis 1979, Stuart Margolis and John Rhodes 2022-2023.}

\item{{\bf $\GM$ Semigroups Built From Character Tables of Abelian Groups}\label{chartab}

In this section we build a $\GM$ semigroup $S(2)$ of complexity 1 using the character theory of a cyclic group of order $4$. We show that $S(2)$ has complexity 1 by showing that $\RLM(S)c=1$ and constructing a flow over $RZ(4)^1$. We indicate how the construction can be generalized to a semigroup $S_{n}$ by using the character theory of a cyclic group of order $2^n$.

We write the cyclic group of order $n$ with generator $x$ multiplicatively. That is, $Z_{n}=\{x^{i} \mid i = 0,\ldots , n-1\}$.
 We define the ``abstract'' character table of $Z_{n}$ to be the matrix $C_n$ with rows and columns indexed by $\{0,\ldots , n-1\}$ and such that 
$C_{n}(k,l)=x^{kl}$. If we substitute the complex $n^{th}$ primitive root of unity, $e^{2\pi i/n}$ for $x$ 
 we obtain the complex character table $\chi_{n}$ of $Z_{n}$. Thus in $C_n$ we are effectively identifying $Z_{n}$ with its complex character group. 

We define $Ch_{n}$ to be the simple semigroup $\mathcal{M}(Z_{n},\{0,\ldots , n-1\},\{0,\ldots , n-1\},C_{n})$. Since the rows and columns of $Ch_n$
are not proportional, $Ch_{n}$ is a $\GM$ simple semigroup with maximal subgroup $Z_n$. By Theorem 7.1 of \cite{FlowsI}, $Ch_n$ is a subsemigroup of the wreath product $Z_{n} \wr (n,RZ(n))$. 

Let $X$ be the $n \times n$ permutation matrix with indices 
$\{0,\ldots , n-1\}$ such that $X(i,i+1)=1$, where we take $i$ modulo $n$. That is, $X$ is the permutation matrix of the cycle $(01\ldots n-1)$ acting on the right. Let $Y$ be the $n \times n$ diagonal monomial matrix with $Y(i,i) = x^{i}$. Direct matrix multiplication shows that
$XC_{n} = C_{n}Y$. Effectively, shifting the rows up 1 is the same as multiplying column $i$ by $x^{i}$. Therefore, the pair $(X,Y)$ is in the translational hull of $Ch_{n}$ and we can form the monoid $R_{n} = Ch_{n} \cup Z_{n}$ where the generator $x$ of $Z_{n}$ acts via $X$ on the right of $Ch_{n}$ and by $Y$ on the left. The monoid $Res_{n}=R_{n}^{op}$ is then a submonoid
of $Z_{n} \wr RZ(n)^{1}$ by Theorem Theorem 7.1 of \cite{FlowsI}. We need the opposite monoid because we want a diagonal left action of $x$ on $S(Z_{n},n)$ since SPCs are proportional on the left. We will see below that $Res_{n}$ is the resolution semigroup of the monoid $S_{n}$ that we now define.

We remark that there is a natural representation theoretic meaning to this construction. If we replace $C_{n}$ by the complex character table 
$\chi_{n}$ as explained above, then we can rewrite the equation $XC_{n}=C_{n}Y$ as $Y = C_{n}^{-1}XC_{n}$, where we also replace $x \in Y$ by $e^{2\pi i/n}$. This reflects the well-known fact that the eigenvectors of $\chi_{n}$ give a basis that diagonalizes the right regular representation of $Z_n$ as a direct sum of one copy of each irreducible representation.

We now define a semigroup $S_{2}$ with the following properties:

\begin{enumerate}

\item{$S_{2}$ is a $\GM$ semigroup whose 0-minimal $\mathcal{J}$-class has maximal subgroup $Z_{4}$.}

\item{$S_{2}$ has a flow over the semigroup $RZ(4)^{1}$.}

\item{$\RLM(S_{2})c = 1$ and thus $S_{2}c = 1$.}


\end{enumerate} 

The semigroup $S(1)$ is Example \ref{CBIRIP}. 


Let $M_{8}^{T}$ be the $8 \times 16$ matrix who in block form is the matrix
$\begin{bmatrix} \Gamma_{8} & | & I_{8} \end{bmatrix}$ where $I_{8}$ is the identity matrix of size 8. This is the transpose of the matrix $M_{k}$ from Section \ref{BIRIP} in the case $k=8$. 

By Lemma 4.1 of \cite{cremona}, the translational hull $\Omega(\mathcal{M}^{0}(\{1\},A,B,M_{8}^{T})$ of the 0-simple semigroup over the
 trivial semigroup with structure matrix $M_{8}^{T}$ consists of all partial functions $f$ on the vertices of the 8-cycle such that the
 inverse image of a vertex or edge in the image of $f$ is also either a vertex or an edge. This makes it easy to check if such a function
 belongs to the translational hull of this 0-simple semigroup. More generally, if we let $\mathcal{M}_{8}(G)$ be the 0-simple semigroup with
 maximal subgroup $G$ and structure matrix $M_{8}^{T}$, Lemma 4.1 of \cite{cremona} gives a necessary and sufficient condition for a row
 monomial matrix over $G$ to be in the translational hull $\Omega(\mathcal{M}_{8}(G))$. We use this criterion without mentioning it in the
 definition of $S(2)$. 

We now define $S(2)$ as a subsemigroup of the translational hull of $\Omega(\mathcal{M}_{8}(Z_{4}))$. The generators are all the elements of the 0-minimal ideal $\mathcal{M}_{8}(Z_{4})$ together with the following elements:

\begin{enumerate}

\item{Let $a$ be the permutation $(1357)(2468)$ corresponding to adding 2 modulo 8. By our convention, this is shorthand for the row-monomial permutation matrix corresponding to $a$.}

\item{Define $b$ = $8 \rightarrow 1$, $7 \rightarrow 2$. Notice that $a$ maps an element in its domain from an orbit of $a$ to an element in the other orbit.}

\item{We now use the abstract character table $C_4 = \begin{bmatrix} 1 & 1 & 1 & 1 \\ 1&x&x^{2}&x^{3} \\ 1&x^{2}&1&x^{2} \\ 1&x^{3}&x^{2}&x \end{bmatrix}$ to define 3 more generators:
\begin{enumerate}

\item{$s_{x}= 1 \rightarrow 2, 2 \rightarrow x4, 3 \rightarrow x^{2}6, 4 \rightarrow x^{3}8$. That is, $s_{x}$ weights the function that sends 
$1\rightarrow 2, 2\rightarrow 4, 3\rightarrow 6, 4\rightarrow 8$ by the character corresponding to $x$.}

\item{$s_{x^{2}}$ weights $1\rightarrow 2, 2\rightarrow 4, 3\rightarrow 6, 4\rightarrow 8$ by the character corresponding to $x^2$. Thus, 
$s_{x^{2}}= 1 \rightarrow 2, 2 \rightarrow x^{2}4, 3 \rightarrow 6, 4 \rightarrow x^{2}8$.}

\item{$s_{x^{3}}$ weights $1\rightarrow 2, 2\rightarrow 4, 3\rightarrow 6, 4\rightarrow 8$ by the character corresponding to $x^3  $. Thus, 
$s_{x^{3}}= 1 \rightarrow 2, 2 \rightarrow x^{3}4, 3 \rightarrow x^{2}6, 4 \rightarrow x8$.}

\end{enumerate}}

\item{We now define $S_2$ to be the semigroup generated by the 0-simple semigroup $\mathcal{M}_{8}^{T}(Z_4)$, where we take $x$ as the generator of the distinguished copy of $Z_4$ in this semigroup and the elements $a,b,s_{x},s_{x^{2}},s_{x^{3}}$ defined above.}

\end{enumerate}

We first gather some facts about $S_{2}$.

\begin{enumerate}

\item{$S_{2}$ is a subsemigroup of the translational hull of $\mathcal{M}_{8}(Z_{4})$. As mentioned above, this is easily justified by using Lemma 4.1 of \cite{cremona} and Lemma 4.1 of \cite{cremona}. It follows that $S_{2}$ is a $\GM$ semigroup. Furthermore, 
$I(S_{2}) \approx \mathcal{M}^{0}(Z_{4},\{1,\ldots 16\}, \{1 \ldots 8\}, M_{8}^{T})$.}

\item{$S_{2}$ does not have a one-point flow. By Theorem 5.2 of \cite{FlowsI} we must show that $I(S_{2}) \cap (S_{2})_{II}$ is not aperiodic, which we now do. Firstly, by Graham's Theorem \cite{Graham} it follows that $IG(S_{2}) \cap I(S_{2}) \approx \mathcal{M}^{0}(\{1\},\{1,\ldots 16\}, \{1 \ldots 8\}, M_{8}^{T})$ is aperiodic, where $IG(S_{2})$ is the idempotent generated subsemigroup of $S_{2}$. 

Now note that $(4,x^{3},2)s_{x}(4,x^{3},2)=(4,x^{3},2)$. Therefore $(4,x^{3},2)ts_{x} \in (S_{2})_{II}$ for any $t \in (S_{2})_{II}$. Let
$t=(1,1,1) = t^{2} \in  (S_{2})_{II}$. Thus $(4,x^{3},2)ts_{x}=(4,x^{3},2)(1,1,1)s_{x}=(4,x^{3},2) \in (S_{2})_{II}$. By the above, 
$(2,1,4) \in IG(S_{2}) \subseteq (S_{2})_{II}$ and it follows that the non-trivial group element $(4,x^{3},2)(2,1,4)=(4,x^{3},4)$ belongs to $(S_{2})_{II}$. Therefore, $S_2$ does not have a one-point flow.}

\item{We show that $\RLM(S_{2})c=1$. Note that the 0-minimal ideal $I(S_{2})$ of $\RLM(S_{2})$ is the aperiodic 0-simple semigroup with structure matrix $M_{8}^{T}$. Therefore $\RLM(S_{2})c = (\RLM(S_{2})/I(S_{2}))c$ by the fundamental lemma of complexity \cite[Chapter 4]{qtheory}. Since every element in this Rees quotient acts by partial 1-1 maps on $I(S_{2})$, we have that $\RLM(S_{2})/I(S_{2})$ is a subsemigroup of the Symmetric Inverse Semigroup on $B$. Since every inverse semigroup has complexity at most 1, we are done, given that the element $a$ generates a non-trivial group.}

\end{enumerate}
 
We now compute some states in the evaluation transformation semigroup of $S_{2}$. 

$(\{ 1\}/\langle 1\rangle )a^{\omega+*}=
\{1 \  \mid 3 \ \mid 5 \ \mid 7\}/\langle 1\ \mid \ 1 \mid \ 1 \mid \ 1\rangle $, $\{(1 \  \mid 3 \ \mid 5 \  \mid 7\}/\langle 1\ \mid \ 1 \mid \ 1 \mid \ 1\rangle )b = \{2\}/\langle 2 \rangle $, $(\{2\}/\langle 2 \rangle )a^{\omega+*}=
\{2 \ \mid 4 \ \mid 6 \ \mid 8\}/\langle 1\ \mid \ 1 \mid \ 1 \mid \ 1\rangle $, $(\{2 \ \mid 4 \ \mid 6 \ \mid 8\}/\langle 1\ \mid \ 1 \mid \ 1 \mid  \ 1\rangle )b = \{1\}/ \langle 1\rangle $.

It follows that if we let $D=(b(a^{\omega+*}))^{\omega+*}$, then $(\{1\}/\langle 1\rangle )aD=\{1 \ 2 \ 3 \ 4 \ 5 \ 6 \ 7 \ 8\}/\rangle 1\ 1\ 1\ 1\ 1\ 1\ 1\ 1 \rangle $, where we have silently used Lemma 4.14.29 of \cite{qtheory} known as the Tie-Your-Shoes Lemma to show that we have one block. Therefore 
$\sigma_{4}=\{1 \ 2 \ 3 \ 4 \ 5 \ 6 \ 7 \ 8/\langle 1\ 1\ 1\ 1\ 1\ 1\ 1\ 1\rangle $ is a State as are:

\begin{itemize}

\item{$\sigma_{1}=(\{1 \ 2 \ 3 \ 4 \ 5 \ 6 \ 7 \ 8/\langle 1\ 1\ 1\ 1\ 1\ 1\ 1\ 1\rangle s_{x}=\{2 \ 4 \ 6 \ 8\}/\langle 1 \ x\ x^{2}\ x^{3}\rangle )$}

\item{$ \sigma_{2}=(\{1 \ 2 \ 3 \ 4 \ 5 \ 6 \ 7 \ 8\}/\langle 1\ 1\ 1\ 1\ 1\ 1\ 1\ 1\rangle )s_{x^{2}}=\{2 \ 4 \ 6\ 8\}/\langle 1\ x^{2}\ 1\ x^{2}\rangle$ and}

\item{$\sigma_{3}=(\{1 2 3 4 5 6 7 8\}/\langle 1\ 1\ 1\ 1\ 1\ 1\ 1\ 1\rangle )s_{x^{3}}=\{2 \ 4 \ 6\ 8\}/\langle 1\ x^{3}\ x^{2}\ x\rangle $.}
\end{itemize}

We now define a flow $F:(\{1,2,3,4\},RZ(4)^{1}) \rightarrow Rh_{B}(S_{2})$. We let $F(i)=\sigma_{i}, i =1,2,3,4$. We cover $a$ by the identity of $RZ(4)^{1}$, $b$ by $\bar{4}$, the constant map to $4$ and $s_{x^{i}}$ by $\bar{i}, i=1,2,3$. Finally we cover any element in $I(S_{2})$ by $\bar{4}$.

We verify that $F$ is a flow. We first claim that for $i=1,2,3,4$, $\sigma_{i}a=\sigma_{i}$. This is clear for $i=4$, since $a$ permutes $1, \ldots ,8$ and acts as the identity on the $Z_4$ coordinate. We have $\sigma_{1}a=(\{2\ 4\ 6\ 8\}/\langle 1\ x\ x^{2}\ x^{3}\rangle)a=\{2\ 4\ 6\ 8\}/\langle x^{3}\ 1\ x\ x^{2}\rangle $. But the cross-section $\langle x^{3}\ 1\ x\ x^{2}\rangle  =\langle 1\ x\ x^{2}\ x^{3}\rangle $ since $x(x^{3}\ 1\ x\ x^{2})=(1\ x\ x^{2}\ x^{3})$. Therefore $F_{1}a = \sigma_{1}a=\sigma=F_{1}$. A similar calculation shows that $F_{2}a=F_{2}$ and $F_{3}a=F_{3}$. This is in fact using the elementary character theory of $Z_{4}$ as noted above. Therefore the conditions for a flow are verified for $a$. For $i=1,2,3$, 
$F(4)\sigma_{x^{i}}=(\{1 2 3 4 5 6 7 8\}/\langle 1\ 1\ 1\ 1\ 1\ 1\ 1\ 1\rangle )s_{x^{i}}=\sigma_{i}=F_{i}$ and is undefined on $F_{1},F_{2}$ and $F_{3}$. The image of of  $(a,g,b) \in I((S))$ on $F_{i}$ is either empty or less than or equal to $F_{4}$. This completes the proof that $F$ is a flow.

We therefore have that $(Z_{4} \times \{1,2,3,4\},S_{2}) \prec (Z_{4} \wr Sym(4) \wr (\{1,2,3,4\},RZ(4)^{1}) \times \RLM(S_{2})$ by Theorem \ref{Apuniform} and thus $S_{2}c=1$.

The definition of  $S_2$ can be generalized to a semigroup $S_n$ using the character table of $Z_(2^{n})$ for all $n>0$. See \cite{FlowsI} Section 8.

History: Stuart Margolis and John Rhodes, 2025 \cite{FlowsI}.
}
}
\end{enumerate}
\appendix

\section{APPENDIX: Semilocal Theory, Group Mapping and Right Letter Mapping Semigroups}\label{Semiloc}

Green's relations and Rees Theorem give the local structure of semigroups. Global semigroup theory is involved with decomposition theorems like the Krohn-Rhodes theorem and related tools to understand them. Semilocal theory is involved with how a semigroup acts on its Green classes. The main example is the Schutzenberger representation. We review the theory here. See Chapters 7 and 8 of \cite{Arbib} and 
Chapter 4 of \cite{qtheory} for more detailed information.
\bd

\begin{enumerate}
    
\item{A semigroup $S$ is right (left) transitive if $S$ has a faithful transitive right (left) action on some set $X$.}

\item{A semigroup $S$ is bi-transitive if $S$ has a faithful transitive representation on the right of some set $X$ and $S$ has a faithful transitive representation on the left of some set $Y$. We emphasize that we do not assume that $X=Y$.}

\end{enumerate}
\ed

\brm

In the literature these are called Right Mapping (denoted $\RM$), Left Mapping ($\LM$) and Generalized Group Mapping ($\GGM$) semigroups respectively. We will use these names in this document.

\erm

The most important example of a right transitive action is the Schutzenberger representation $Sch(R)$ on an $\mathcal{R}$-class $R$ of $S$ whose definition we recall.

\medskip
For $r \in R, s \in S$, define 
    
		$$r\cdot s = \begin{cases}
      rs, & \text{if}\ rs \in R \\
      \text{undefined} & \text{otherwise}
    \end{cases}$$

There is the obvious dual notion of the left Schutzenberger representation of a semigroup on an $\mathcal{L}$-class.

\bl
 $S$ is a $\RM$, ($\LM$, $\GGM$) semigroup if and only if $S$ has a unique $0$-minimal regular ideal $I(S) \approx M^{0}(A,G,B,C)$ such that the right
 (left, right and left) Schutzenberger representation of $S$ is faithful on any $\mathcal{R}( \mathcal{L}, \mathcal{R} {\text{ and }} \mathcal{L})$-class $R (L, R {\text{ and }} L)$ class in $I(S) -\{0\}$.  
\el

\bd

A $\GGM$ semigroup is called {\em Group Mapping} written $\GM$ if the group $G$ in $I(S)$ is non-trivial.

\ed

Let $S$ be a $\GM$ semigroup. We identify a fixed $\mathcal{R}$-class in $I(S)$ with the set $G \times B$. We obtain a faithful transformation semigroup $(G \times B,S)$.

The action of $S$ on $G \times B$ induces an action of $S$ on $B$.  The faithful image of this action is called the 
{\em Right Letter Mapping} (written $\RLM(S)$) image of $S$.  We thus obtain a faithful transformation semigroup $(B,RLM(S))$. We remark that $\RLM(S)$ has a 0-minimal regular ideal that is aperiodic and is the image of $I(S)$ under the morphism from $S$ to $\RLM(S)$. In particular, for any $\GM$ semigroup $S$, we have $|\RLM(S)| < |S|$. This allows for inductive proofs based on cardinality.

\bl

Let $S$ be a $GM$ semigroup. Then the following hold.

\begin{enumerate}

\item{$(G \times B,S)$ is a subsemigroup of $G \wr (B,\RLM(S))$. Therefore, $Sc \leq 1+ \RLM(S)c$. }


\item{For every semigroup $T$, there is a $GM$ quotient semigroup $S$ such that $Tc=Sc$.} 

\end{enumerate}

\el

By induction on cardinality, we are lead to the following reduction theorem.

\bt
The question of decidability of complexity can be reduced to the case that $S$ is a $\GM$ semigroup and whether $Sc=\RLM(S)c$ or 
$Sc = 1+RLM(S)c$.

\et

\brm

It is known that if $S$ is a completely regular semigroup, that is, a semigroup that is a union of its subgroups, then if $S$ is $\GM$ then it is always true that $Sc = 1+\RLM(S)c$ and there are completely regular semigroups of arbitrary complexity \cite[Chapter 9]{Arbib}. On the other hand, if $S$ is an inverse $\GM$ semigroup such that $\RLM(S)$ is not aperiodic, then $Sc=\RLM(S)c = 1$. The semigroup of all row and column monomial matrices over a non-trivial group $G$ is example of such an inverse semigroup.

\erm

\section{APPENDIX: Connections Between the Set-Partition Lattice and the Rhodes Lattice} \label{SetPart}

Let $S$ be a $\GM$ semigroup with $0$-minimal ideal ${\cal M}^0(G, A, B, C)$. The definition of flow gives a map from the state set 
of an automaton to a lattice associated with $S$. In the literature, there are two such lattices. We define these lattices and give the connections 
between them and show that they give equivalent notions of flow. The first is the set-partition lattice $\operatorname{SP}(G \times B)$. This is the lattice 
whose elements are all pairs $(Y,\Pi)$, where $Y$ is a subset of $G\times B$ and $\Pi$ is a partition on $Y$. Here $(Y,\Pi) \leqslant (Z,\Theta)$ if 
$Y \subseteq Z$ and for all $y \in Y$, the $\Pi$ class of $y$ is contained in the $\Theta$ class of $y$. 

The second lattice is the Rhodes lattice $\operatorname{Rh}_{B}(G)$. We review the basics. For more details see~\cite{AmigoDowling}. 
Let $G$ be a finite group and $B$ a finite set. A partial partition on $B$ is a partition $\Pi$ on a subset $I$ of $B$. We also consider the 
collection of all functions $F(B,G)$, $f \colon I \rightarrow G$ from subsets $I$ of $B$ to $G$. The group $G$ acts on the left of $F(B,G)$ by $(gf)(b) = gf(b)$
for $f \in F(B,G), g \in G, b \in \operatorname{Dom}(f)$. An element $\{gf \mid f \colon I \rightarrow G, I \subseteq X, g \in G\}$ of the quotient set 
$F(B,G)/G$ is called a cross-section 
with domain $I$. It should be thought of as a projectification of a cross-section of the projection from $I$ to $B$ in the usual topological sense. An
SPC (Subset, Partition, Cross-section) over $G$ is a triple $(I,\Pi,f)$, where $I$ is a subset of $B$, $\Pi$ is a partition of $I$ and $f$ is a collection of 
cross-sections one for each $\Pi$-class $\pi$ with domain $\pi$. If the classes of $\Pi$ are $\{\pi_{1},\pi_{2}, \ldots, \pi_{k}\}$, then we sometimes write 
$\{(\pi_{1},f_{1}), \ldots, (\pi_{k},f_{k})\}$, where $f_{i}$ is the cross-section associated to $\pi_{i}$. For brevity we denote this set of cross-sections by 
$[f]_{\pi}$.  We let $\operatorname{Rh}_{B}(G)$ denote the set of all SPCs on $B$ over the group $G$ union a new element 
$\Longrightarrow\Longleftarrow$  that we call {\em contradiction} and is the top element of the lattice structure on $\operatorname{Rh}_{B}(G)$.
Contradiction occurs  because the join of two SPCs need not exist. In this case we say that the contradiction is their join.

The partial order on $\operatorname{Rh}_{B}(G)$ is defined as follows. We have $(I,\pi,[f]_{\pi}) \leqslant (J,\tau,[h]_{\tau})$ if:
\begin{enumerate}
\item
{$I \subseteq J$;} 
\item
{Every block of $\pi$ is contained in a (necessarily unique) block of $\tau$;} 
\item{if the $\pi$-class $\pi_{i}$ is a subset of the $\tau$-class $\tau_{j}$, then the restriction of $h$ to $\pi_{i}$ equals $f$ restricted to $\pi_{i}$ as elements of $F(B,G)/G$. That is, $[h|_{\pi_{i}}] = [f|_{\pi_{i}}]$.}
\end{enumerate}

See \cite[Section 3]{AmigoDowling} for the definition of the lattice structure on $\operatorname{Rh}_{B}(G)$. The underlying set of the Rhodes lattice 
$\operatorname{Rh}_{B}(G)$ minus the contradiction is the set underlying the Dowling lattice on the same set and group. The Dowling lattice 
has a different partial order. For the connection between Rhodes lattices and Dowling lattices see \cite{AmigoDowling}.

We note that $\operatorname{SP}(G \times B)$ is isomorphic to the Rhodes lattice $\operatorname{Rh}_{G\times B}(1)$ of the trivial group over the 
set $G \times B$. We need only note that a cross-section to the trivial group is a partial constant function to the identity and can be omitted, leaving 
us with a set-partition pair. There are no contradictions for Rhodes lattices over the trivial group, and in this case the top element is the pair 
$(G \times B,(G \times B)^{2})$. Despite this, we prefer to use the notation $\operatorname{SP}(G \times B)$ instead of $\operatorname{Rh}_{G\times B}(1)$.

Conversely, we can find a copy of the meet-semilattice of $\operatorname{Rh}_{B}(G)$ as a meet subsemilattice of $\operatorname{SP}(G\times B)$. 
We begin with the following important definition.

\bd
A subset $X$ of $G \times B$ is a cross-section if whenever $(g,b),(h,b) \in X$, then g = h. That is, $X$ defines a cross-section of the projection 
$\theta \colon X \rightarrow B$. Equivalently, $X^{\rho} \subseteq B \times G$, the reverse of $X$, is the graph of a partial function 
$f_{X} \colon B \rightarrow G$.  An  element $(Y,\Pi) \in \operatorname{SP}(G \times B)$ is a cross-section if  every partition class $\pi$ of $\Pi$ 
is a cross-section.  
\ed

From the semigroup point of view, a cross-section is a partial transversal of the $\mathcal{H}$-classes of the distinguished $\mathcal{R}$-class, 
$R = G \times B$ of a $\GM$ semigroup. That is, the $\mathcal{H}$-classes of $R$ are indexed by $B$ and a cross-section picks at most one 
element from each $\mathcal{H}$-class. 

An element $(Y,\Pi)$ that is not a cross-section is called a {\em contradiction}. That is, $(Y,\Pi)$ is a contradiction if some $\Pi$-class $\pi$ contains 
two elements $(g,b),(h,b)$ with $g \neq h$. We note that the set of cross-sections is a meet subsemilattice of $\operatorname{SP}(G \times B)$ and 
the set of contradictions is a join subsemilattice of $\operatorname{SP}(G \times B)$.

We will identify $B$ with the subset $\{(1,b)|b \in B\}$ of $G \times B$, which as above we think of as a system of representatives of the 
$\mathcal{H}$-classes of the distinguished $\mathcal{R}$-class. Note then that $G\times B$ is the free left $G$-act on 
$B$ under the action $g(h,b)=(gh,b)$. This action extends to subsets and partitions of $G\times B$ . Thus $\operatorname{SP}(G\times B)$ 
is a left $G$-act. An element $(Y,\Pi)$ is {\em invariant} if $G(Y,\Pi) = (Y,\Pi)$. It is easy to see that $(Y,\Pi)$ is invariant if and only if:

\begin{enumerate}
 \item{$Y=G \times B'$ for some subset $B'$ of $B$.}  

 \item{For each $\Pi$-class $\pi$, $G\pi \subseteq \Pi$ and is a partition of $G\times B''$ where $B''\subseteq B'$.}
\end{enumerate}

Thus $(Y,\Pi)$ is invariant if and only if $Y=G\times B'$ for some subset $B'$ of $B$ and there is a partition $B_{1},\ldots B_{n}$ of $B'$ such that for all $\Pi$ classes $\pi$, $G\pi$ is a partition of $G\times B_{i}$ for a unique $1 \leqslant i \leq n$. 

Let $\operatorname{CS}(G\times B)$ be the set of invariant cross-sections $(Y,\Pi)$ in $\operatorname{SP}(G \times B)$. $\operatorname{CS}(G\times B)$ is a 
meet-subsemilattice of $\operatorname{SP}(G\times B)$. We claim it is isomorphic to the meet-semilattice of $\operatorname{Rh}_{B}(G)$. Indeed, 
let $(G\times B',\Pi)$ be an invariant cross-section. Using the notation above, pick  $\Pi$-classes, $\pi_{1}, \ldots , \pi_{n}$ such that $G\pi_{i}$ is a 
partition of $G\times B_{i}$ for $i=1, \ldots n$. Since $\pi_{i}$ is a cross-section its reverse is the graph of a function $f_{i} \colon B_{i}\rightarrow G$. We map 
$(G\times B',\Pi)$ to $(B',\{B_{1},\ldots , B_{n}\},[f]_{\{B_{1},\ldots , B_{n}\}}) \in \operatorname{Rh}_{B}(G)$, where the component of $f$ on $B_{i}$ is 
$f_i$. It is clear that this does not depend on the representatives $\pi$ and we have a well-defined function from 
$F \colon \operatorname{CS}(G \times B) \rightarrow \operatorname{Rh}_{B}(G)$. It is easy to see that $F$ is a morphism of meet-semilattices.

Conversely, let $(B',\Theta,[f]) \in \operatorname{Rh}_{B}(G)$. We map $(B',\Theta,[f])$ to $(G \times B',\Pi) \in \operatorname{CS}(G \times B)$,
where $\Pi$ is defined as follows. If $\theta$ is a $\Theta$-class, let $\widehat{\theta} =\{(bf,b) \mid b \in \theta\}$. Let $\Pi$ be the collection of all subsets 
of the form $g\widehat{\theta}$ where $\theta$ is a $\Theta$-class and $g \in G$. The proof that $(G \times B',\Pi) \in \operatorname{CS}(G \times B)$ 
and that this assignment is the inverse to $F$ above and gives an isomorphism between the meet-semilattice of $\operatorname{Rh}_{B}(G)$ and 
$\operatorname{CS}(G \times B)$ is straightforward and left to the reader. Furthermore, if the join of two SPC in $\operatorname{Rh}_{B}(G)$ 
is an SPC (that is it is not the contradiction) then this assignment preserves joins. Therefore we can identify the lattice $\operatorname{Rh}_{B}(G)$ 
with $\operatorname{CS}(G \times B)$ with a new element $\Longrightarrow\Longleftarrow$ added when we define the join of two elements of 
$\operatorname{CS}(G \times B)$ to be $\Longrightarrow\Longleftarrow$ if their join in $\operatorname{SP}(G\times B)$ is a contradiction. 
We use $\operatorname{CS}(G\times B)$ for this lattice as well.

We record this discussion in the following proposition.

\bp\label{Updown}
The set-partition lattice $\operatorname{SP}(G\times B)$ is isomorphic to the Rhodes lattice $\operatorname{Rh}_{G\times B}(1)$. 
The Rhodes lattice $\operatorname{Rh}_{B}(G)$ is isomorphic to the lattice $\operatorname{CS}(G\times B)$.
\ep

\section{APPENDIX: Flows and the Flow Decomposition Theorem}\label{Flows}

We review the definition of flows from  an automaton with alphabet $X$  to the set-partition $SP(G\times B)$ and to the Rhodes lattice $Rh_{B}(G)$. For more details, see \cite[Sections 2-3]{Trans}. We first need to recall details about the Type II subsemigroup of a semigroup.

\subsection{The Type II Subsemigroup and the Tilson Congruence}\label{redux.sec}

In this subsection we review the type II subsemigroup $S_{II}$ of a semigroup. It plays an important role in finite semigroup theory and a central role in flow theory. It was first defined in \cite{lowerbounds2}. In that paper it was proved that it is decidable if a regular element of a semigroup $S$ belongs to $S_{II}$. In particular, it followed that if $S$ is a regular semigroup, then membership in $S_{II}$ is decidable. Later 
Ash \cite{Ash} proved that membership in $S_{II}$ is decidable for all finite semigroups $S$.

The type II subsemigroup $S_{II}$ of $S$ is the smallest subsemigroup of $S$ containing all idempotents and closed under weak conjugation: if $xyx = x$ for $x,y \in S$, then $xS_{II}y \cup yS_{II}x \subseteq S_{II}$. Membership in $S_{II}$ is clearly decidable from this definition. Its importance stems from the following Theorem, where we give its original, not a priori decidable definition.

\bt

Let $S$ be a finite semigroup. Then $S_{II}$ is the intersection of all subsemigroups of the form $1\phi^{-1}$ where $\phi:S \rightarrow G$
is a relational morphism from $S$ to a group $G$.

\et

The aforementioned decidability results in \cite{lowerbounds2, Ash} prove that the two definitions we give define the same subsemigroup of a semigroup $S$. Here are some important properties of the type II subsemigroup. See \cite{qtheory} for proofs.

\bt\label{typeIIprop}

\begin{enumerate}

\item{Let $f:S \rightarrow T$ be a morphism between semigroups $S$ and $T$. If $s \in S_{II}$, then $sf \in T_{II}$. Therefore the restriction of $f$ to $S_{II}$ defines a morphism $f_{II}:S_{II} \rightarrow T_{II}$. This defines a functor on the category of finite semigroups.}

\item{If $S$ divides $T$, then $S_{II}$ divides $T_{II}$.}

\item{Assume that a semigroup $S$ divides a semidirect product $T*G$ where $T$ is a semigroup and $G$ is a group. Then $S_{II}$ divides $T$.}

\end{enumerate}

\et





A {\em congruence} on a transformation semigroup $(Q,S)$ is an equivalence relation $\approx$ on $Q$ such that if $q\approx q'$ and for $s \in S$, and both $qs$ and $q's$ are defined, then $qs\approx q's$. Every $s \in S$ defines a partial function on $\faktor{Q}{\approx}$ by 
$[q]_{\approx}s=[q's]_{\approx}$ if  $q's$ is defined for some $q' \in [q]_{\approx}$. The quotient $\faktor{(Q,S)}{\approx}$ has states $\faktor{Q}{\approx}$ and semigroup $T$ the semigroup generated by the action of all $s \in S$ on $\faktor{Q}{\approx}$. We remark that $T$ is not necessarily a quotient semigroup of $S$, but is in the case that $(Q,S)$ is a transformation semigroup of total functions. 

A congruence $\approx$ is called {\em injective} if every $s \in S$ defines a partial 1-1 function on $\faktor{Q}{\approx}$. It is easy to see that the intersection of injective congruences is injective. Therefore, there is a unique minimal injective congruence $\tau$ on any transformation semigroup $(Q,S)$. We call $\tau$ the Tilson congruence on a transformation semigroup because of the following proved in \cite{Redux}. It is central to the theory of flows.

\bt\label{redux}

Let $(G \times B,S)$ be a $\GM$ transformation semigroup. Then the minimal injective congruence $\tau$ on $S$ is defined as follows: 
$(g,b) \tau (g',b')$ if and only if there are elements $s,t \in S_{II}$ such that $(g,b)s=(g',b')$ and $(g',b')t=(g,b)$.

\et
 
\brm

\item{We can state this theorem by $(g,b)S_{II}=(g',b')S_{II}$. That is, $(g,b)$ and $(g',b')$ define the same ``right coset'' of $S_{II}$ on $G \times B$.}

\item{The proof in \cite{Redux} works on an arbitrary transitive transformation semigroup. We stated it in the case of $\GM$ transformation semigroups because that's how we will use it in this document.}

\item{Theorem \ref{redux} can be used to greatly simplify the proof in \cite{lowerbounds2} for decidability of membership in $S_{II}$ for regular elements of an arbitrary  semigroup.}

\erm

For later use we record the following corollary to Theorem \ref{redux}.

\bc

Let $(G \times B,S)$ be a $\GM$ transformation semigroup. Then $(g,b) \tau (g',b')$ if and only if there are elements 
$s,t \in S_{II} \cap I(S)$ such that $(g,b)s=(g',b')$ and $(g',b')t=(g,b)$. That is, we can choose the elements $s$ and $t$ in Theorem \ref{redux} to be in the 0-minimal ideal of $S$.

\ec

\proof The condition is sufficient by Theorem \ref{redux}. Conversely, assume that $(g,b)\tau (g',b')$. Then there are elements $s,t \in S_{II}$ such that $(g,b)s=(g',b')$ and $(g',b')t=(g,b)$. Since $I(S)$ is a 0-simple semigroup, there are idempotents $e,f \in I(S)$ such that 
$(g,b)e=(g,b)$ and $(g',b')f=(g',b')$. Therefore, $(g,b)es=(g,b)s=(g',b')$ and similarly $(g',b')ft =(g,b)$. Since $e,f$ are idempotents we have $es,ft \in S_{II}$. As $I(S)$ is a 0-minimal ideal, we also have $es,ft \in I(S)$. \Qed

%
%

\subsection{Definition of Flows and Their Properties}
Let $(G \times B,S)$ be the transformation semigroup associated to a $\GM$ semigroup $S$ and let $X$ be a generating set for $S$. By a deterministic automaton we mean an automaton such that each letter defines a partial function on the state set.

\bd

Let $\mathcal{A}$ be a deterministic automaton with state set $Q$ and alphabet $X$. 
A {\em flow} to the lattice $SP(G\times B)$ on $\mathcal{A}$ is a function $f:Q \rightarrow SP(G \times B)$ such that for each $q \in Q, x \in X$, with $qf =(Y,\pi)$ and $(qx)f = (Z,\theta)$ we have:

\begin{enumerate}

\item{For all $(g,b) \in G \times B$, there is a $q \in Q$ such that $(g,b) \in qf$.}
 
\item{$Yx \subseteq Z$.}

\item{Multiplication by $x$ considered as an element of $S$ induces a partial 1-1 map from $Y/\pi$ to $Z/\theta$.

That is, for all $y,y' \in Y$ we have $y,y'$ are in a $\pi$-class if and only if $yx,y'x$ are in a $\theta$-class whenever $yx,y'x$ are both defined.}

\item{{\bf The Cross Section Condition:} For all $q \in Q$, $qf$ is a cross-section. That is, for all $g,h \in G, b \in B$ if $(g,b),(h,b) \in qf$ it follows that $g=h$.}

\end{enumerate}
\ed

We use Proposition \ref{Updown} to give a definition of a flow to the Rhodes lattice $Rh_{B}(G)$.

\bd

Let $\mathcal{A}$ be a deterministic automaton with state set $Q$ and alphabet $X$. 
A {\em flow} to the lattice $Rh_{B}(G)$ is a flow to $SP(G\times B)$ such that for each state $q$, $qf \in CS(G\times B)$. That is, $qf$ is a $G$-invariant cross-section.
   
\ed

\brm

    It follows from the original statement of the Presentation Lemma \cite{AHNR.1995} that if $S$ has a flow with respect to some automaton over $SP(G\times B)$, then it has a flow from the same automaton over $Rh_{B}(G)$. The proof of the equivalence of the Presentation Lemma and Flows in Section 3 of \cite{Trans} works as well for the Presentation Lemma in the sense of \cite{AHNR.1995}. 

\erm

We will use all the terminology and concepts from \cite{Trans}. If $\mathcal{A}$ is an automaton with alphabet $X$ and state set $Q$, recall that its completion is the automaton $\mathcal{A}^{\square}$ that adds a sink state $\square$ to $Q$ and declares that $qx = \square$ if $qx$ is not defined in $\mathcal{A}$. A flow on $\mathcal{A}$ is a complete flow if on $\mathcal{A}^{\square}$ extending $f$ by letting $\square f = (\emptyset, \emptyset)$ the bottom of both the set-partition and Rhodes lattices remains a flow and furthermore, for all $(g,b) \in G \times B$, there is a $q\in Q$ such that $(g,b) \in Y$, where $qf=(Y,\Pi)$. All flows in this paper will be complete flows.

 Flows are related to the Presentation Lemma \cite{AHNR.1995}, \cite[Section 4.14]{qtheory}. They give a necessary and sufficient condition for $Sc =\operatorname{RLM}(S)c$ where $S$ is a $\operatorname{GM}$ semigroup.  See Section 3 of \cite{Trans} for a proof of the following Theorem

\bt\label{PLflow}[The Presentation Lemma-Flow Version]

Let $(G \times B,S)$ be a $GM$ transformation semigroup with $S$ generated by $X$. Let $k> 0$ and assume that $\operatorname{RLM}(S)c = k$.  Then $Sc = k$ if and only if there is an $X$ automaton $\mathcal{A}$ whose transformation semigroup $T$ has complexity strictly less than $k$ and a complete flow $f:Q \rightarrow Rh_{B}(G)$.


\et

 Flows were preceded by the Presentation Lemma \cite{AHNR.1995}, \cite[Section 4.14]{qtheory}. The Presentation Lemma was shown to be equivalent to the existence of an appropriate flow in \cite[Section 3]{Trans}. There are three versions of the Presentation Lemma and its relation to Flow Theory in the literature \cite{AHNR.1995}, \cite[Section 4.14]{qtheory}, \cite{Trans}. These have very different formalizations and it is not clear how to pass from one version to another. The terminology is different as well. For example, the definition of cross-section in each of these references, as well as the one we use in this paper is different.  The following Theorem gives a unified approach to these topics. It summarizes known results in the literature and is meant to emphasize the strong connections between slices, flows and the presentation lemma and the corresponding direct product decomposition. \cite{FlowsI} for a proof. We use some of the results that we've proved in previous sections. For background on the derived transformation semigroup and the derived semigroup theorem see \cite{Eilenberg, qtheory}.

\bt \label{uniform}

Let $(G \times B,S)$ be $\GM$ and assume that $\RLM(S)c \leq n$. Then the following are equivalent:

\begin{enumerate}

\item{$Sc \leq n$.}

\item{There is an aperiodic relational morphism $\Theta:S \rightarrow H\wr T$, where $H$ is a group and $Tc \leq n-1$.}

\item{There is a relational morphism $\Phi:S \rightarrow T$ where $Tc \leq n-1$ and such that the Derived Transformation semigroup $D(\Phi)$ is in $Ap*Gp$.}

\item{There is a relational morphism $\Phi:S \rightarrow T$ where $Tc \leq n-1$ such that the Tilson congruence $\tau$ on the Derived Transformation semigroup $D(\Phi)$ is a cross-section.}

\item{$(G \times B,S)$ admits a flow from a transformation semigroup $(Q,T)$ with $(Q,T)c \leq n-1$.}

\item{$S \prec (G \wr (Sym(B)) \wr T) \times \RLM(S)$ for some transformation semigroup $T$ with $Tc \leq n-1$.}

\end{enumerate}

\et

%
%

We emphasize the case for complexity 1 and aperiodic flows as this is used in many of our examples.

\bt \label{Apuniform}

Let $(G \times B,S)$ be $\GM$ and assume that $\RLM(S)c \leq 1$. Then the following are equivalent:

\begin{enumerate}

\item{$Sc=1$.}

\item{There is an aperiodic relational morphism $\Theta:S \rightarrow H\wr T$, where $H$ is a group and $T$ is aperiodic.}

\item{There is a relational morphism $\Phi:S \rightarrow T$ where $T$ is aperiodic such that the Derived transformation semigroup $D(\Phi)$ is in $Ap*Gp$.}

\item{There is a relational morphism $\Phi:S \rightarrow T$ where $T$ is aperiodic such that the Tilson congruence $\tau$ on the Derived transformation semigroup $D(\Phi)$ is a cross-section.}

\item{$(G \times B,S)$ admits an aperiodic flow.}

\item{$S \prec (G \wr (Sym(B)) \wr T) \times \RLM(S)$ for some aperiodic semigroup $T$.}

\end{enumerate}

\et

\section{APPENDIX The Flow Monoid and the Evaluation Transformation Semigroup}\label{Eval}

\subsection{The Monoid of Closure Operations}\label{Closops}

In this Appendix we gather definitions and results from \cite{Trans}. For more details the reader is strongly urged to consult this paper. Since all the computations we do on examples in this document are done in the Evluation Transformation Semigroup defined in Section 5 of \cite{Trans} our goal is to summarize the background material in that paper needed to define this object. We begin with the definition of the monoid of closure operations on the direct product $L^{2} = L \times L$ of a lattice $L$ with itself. 

Let $L$ be a lattice and let $L^{2} = L \times L$. Let $f$ be a closure operator on $L^2$. By definition this means that $f$ is an order preserving, extensive (that is, for all $(l_{1},l_{2}) \in L^{2}, (l_{1},l_{2}) \leqslant (l_{1},l_{2})f$), idempotent function on $L^2$. A {\em stable pair} for $f$ is a closed element of $f$. Thus a stable pair $(l,l') \in L^2$ is an element such that $(l,l')f = (l,l')$.  The stable pairs of $f$ are a meet closed subset of $L \times L$. Conversely, each meet closed subset of $L \times L$ is the set of stable pairs for a unique closure operator on $L \times L$. We will identify $f$ as a binary relation on $L$ whose pairs are precisely the stable pairs. Let $B(L)$ be the monoid of binary relations on the set $L$.  As is well known, $B(L)$ is isomorphic to the monoid $ M_{n}(\mathcal{B})$ of $n \times n$ matrices over the 2-element Boolean algebra $\mathcal{B}$, where $n=|L|$. The Boolean matrix associated to $f$ is of dimension $|L| \times |L|$. It has a 1 in position $(l_{1},l_{2})$ if $(l_{1},l_{2})$ is a stable pair and a 0 otherwise. 

With this identification, the collection $\mathcal{C}(L^{2})$ of all closure operators on $L^2$ is a submonoid of the monoid $B(L)$ of binary relations on $L$~\cite[Proposition 2.5]{Trans}. We thus also consider $\mathcal{C}(L^{2})$ to be a monoid of $|L| \times |L|$ Boolean matrices. Let $L$ be either $Rh_{B}(G)$ or $\operatorname{SP}(G\times B)$. We now recall the definition of some important unary operations on $\mathcal{C}(L^{2})$. 
\bd

Let Let $f \in \mathcal{C}(L^{2})$.

\begin{enumerate}

\item{The domain of $f$ denoted by $\operatorname{Dom}(f)$ is the set $\{x \mid \exists y, (x,y) \in f\}.$}

\item{Define the relation $\overleftarrow{f}$ by $\overleftarrow{f}= \{(x,x) \mid x \in \operatorname{Dom}(f)\}$. $\overleftarrow{f}$ is called {\em back flow along $f$}. See \cite[Remark 2.26]{Trans} for the reason for this terminology.}

\item{The relation $f^{*}$ is defined by $f^{*} = f \cap \{(x,x) \mid x \in X\}$. $f^{*}$ is called {\em the Kleene closure of $R$}. See Section 2 and Section 4 of \cite{Trans} for the reason for this terminology.}

\item{Define the {\em loop of $f$} to be the relation $f^{\omega+*}=f^{\omega}f^{*}$, where $f^{\omega}$ is the unique idempotent in the subsemigroup generated by $f$.}

\end{enumerate}

\ed

At times it is convenient to identify $\overleftarrow{f}$ as $1|_{\operatorname{Dom}(f)}:L \rightarrow L$, the identity function restricted to 
$\operatorname{Dom}(f)$. Similarly, we identify $f^{*}$ as $1|_{\operatorname{Fix{f}}}:L \rightarrow L$, the restriction of the identity to the set of fixed-points of $f$, where $x \in L$ is a fixed point if $(x,x) \in f$. Our use of these will be clear from the context.

\subsection{The 0-Flow Monoid}\label{FMon}

Let $S$ be a $\GM$ semigroup generated by $X$ and let $x \in X$. Let $I(S) = \mathcal{M}^{0}(A,G,B,C)$. We define a binary relation $f_{x}$
on $\operatorname{SP}(G \times B)$ by $((Y,\Pi), (Z,\Theta)) \in f_{x}$ if and only if $Yx \subseteq Z$ and the partial function induced by
right multiplication by $x$, $\cdot x: Y \rightarrow Z$ induces a well-defined partial injective map $\cdot x: Y/\Pi \rightarrow Z/\Theta$. This means that if  $(g,b), (g',b') \in Y$ and $(g,b)x,(g',b')x$ are both defined (and hence in $Z$ ), then $(g,b)\Pi (g',b')$ if and only if $(g,b)x\Theta (g',b')x$. Then $f_{x} \in \mathcal{C}(L^{2})$. See \cite[Proposition 2.22]{Trans}. $f_{x}$ is called the free-flow by $x$.

We now define the 0-flow monoid $M_{0}(L)$ as follows.

\bd\label{Flowops}

Let $L = \operatorname{SP}(G \times B)$. The 0-flow monoid $M_{0}(L)$, is the
smallest subset of $\mathcal{C}(L^{2})$ satisfying the following axioms:

\begin{enumerate}

\item {(Identity) The multiplicative identity $I$ of $\mathcal{C}(L^{2})$ is in $M_{0}(L)$.}

\item{(Points) For all $x \in X$, $f_{x}$ the free-flow along $x$ belongs to $M_{0}(L)$.}

\item{(Products) If $f_{1}, f_{2} \in M_{0}(L)$, then $f_{1}f_{2} \in M_{0}(L)$.}

\item{(Vacuum) If $f \in M_{0}(L)$, then $\overleftarrow{f} \in M_{0}(L)$.}\label{Vac}

\item{(Loops) If $f \in M_{0}(L)$, then $f^{\omega+*}\in M_{0}(L)$.}

\end{enumerate}

\ed

We remark that for each $n \geq 0$, there is an $n$-flow monoid $M_{n}(L)$ defined in \cite{Trans}. The definition above of $M_{0}(L)$ is exactly what is defined in \cite{Trans} and used extensively in \cite{complexity1}. For $n>0$, Axioms (1)-(4) are the same for $M_{n}(L)$ as for $M_{0}(L)$. Axiom (5) restricts the use of the loop operator to $n$-loopable elements \cite[Section 4]{Trans}. In \cite{complexityn}, $n$-loopable elements are replaced by a more restrictive definition and the modified $M_{n}(L)$ plays a crucial role in the main results of \cite{complexityn}.

\subsection{The Evaluation Transformation Semigroup}\label{Ets}

We now defined the Evaluation Transformation Semigroup $\mathcal{E}(L) = (\operatorname{States}, \operatorname{Eval}(L))$ as in \cite{Trans}. Again, because of our interest in this paper on aperiodic flows, we don't define the $n$-Evaluation Transformation Semigroups $E_{n}$ for $n>0$. We begin with the definition of Well-Formed Formualae (WFFs).

\bd

Let $X$ be an alphabet. We define a well-formed formula inductively as follows.

\begin{enumerate}

\item{The empty string $\epsilon$ is a well-formed formula.}

\item{Each letter $x \in X$ is a well-formed formula.}

\item{If $\tau, \sigma$ are well-formed formulae, then
so is $\tau\sigma$.}

\item{If $\tau$ is a well-formed formula that is not a proper power (i.e., not of the form
$\sigma^{n}$ where $n > 1$), then $\tau^{\omega+*}$ is also a well-formed formula.} 

\end{enumerate}

\ed

The set of well-formed
formulae is denoted by $\Omega(X)$. Well-formed formulae will be denoted by Greek
letters. As a convention, if $\tau=\sigma
^{n}$, where $\sigma$ is not a proper power, then we set
$\tau^{\omega+*}=\sigma^{\omega+*}$. In other words, we extract roots before applying the unary operation $\omega+*$.

Let $V = \prod_{f \in M_{0}(L)}\overleftarrow{f}$. $V$ is called the {\em Vacuum}. See \cite{Trans} for an explanation of this terminology. In \cite{Trans}, $V$ is denoted by $\mathcal{F}_{0}$. It is proved in \cite{Trans} that $V$ is an idempotent in $M_{0}(L)$. We now will work exclusively in the subsemigroup $VM_{0}(L)V$ of $M_{0}(L)$. We want to interpret WFFs in $VM_{0}(L)V$.

\bd

Define recursively a partial function $\mathcal{I}: \Omega(X)\rightarrow VM_{0}(L)V$ as follows. 

\begin{enumerate}

\item{$\epsilon\mathcal{I} = V$.}

\item{$x\mathcal{I} = VxV$ for $x \in X$.}

\item{If $\mathcal{I}$ is already defined on $\tau, \sigma \in \Omega(X)$, set $(\tau\sigma)\mathcal{I} = \tau\mathcal{I}\sigma\mathcal{I}$.}

\item{If $\tau \in \Omega(X)$ is not a proper power and $\tau\mathcal{I}$ is defined, set $\tau^{\omega+*}\mathcal{I} = 
(\tau\mathcal{I})^{\omega+*}$.}

\end{enumerate}

\ed

We normally omit $\mathcal{I}$ and assume that a WFF $\tau$  is being evaluated in $VM_{0}(L)V$ according to the definition of $\mathcal{I}$. 
 We first define a new operator on elements of $M_{0}(L)$ called {\em forward flow}. Recall that the bottom of the lattice $\operatorname{SP}(G \times B)$ is the pair $\Box = (\emptyset, \emptyset)$.

\bd

Let $f \in M_{0}(L)$ and let $l \in L$. Let $(l,\Box)f =(l_{1},l_{2})$ We define the forward flow of $f$ denoted by $\overrightarrow{f}$ by 
$l\overrightarrow{f}=l_{2})$. That is, we apply $f$ to $(l,\Box)$ and project to the right-hand coordinate.

\ed

The following is proved in Sections 2 and 4 of \cite{Trans}.

\begin{enumerate}

\item{$\overrightarrow{f}:L \rightarrow L$ is an order preserving function on $L$.}

\item{If $lV=l$, then $(l,\Box)f=(l,l')$ and $l' \in LV$. That is, the left-coordinate of $(l,\Box)f$ is still $l$ and the right-hand coordinate is in $LV$. Therefore $\overrightarrow{f}$ is a well-defined function from $LV$ to $LV$.}

\item{The assignment of $f$ to $\overrightarrow{f}$ defines an action of $VM_{0}(L)V$ on $LV$. It follows that we have a transformation semigroup $(LV,M'_{0}(L))$, where $M'_{0}(L)$ is the image of $VM_{0}(L)V$ on $LV$ under this action.}

\end{enumerate}

We now restrict the action in $(LV,M'_{0}(L))$ to the set of $\operatorname{States}$ defined as follows. Let $(g,b) \in G \times B$. Then 
the element $(\{(g,b)\},\{(g,b)\}) \in L$ is called a {\em point}. In Section 5 of \cite{Trans}, it is proved that every point $p$ satisfies $pV=p$. 

\bd

The set of $\operatorname{States}$ is the smallest subset of $LV$ such that:

\begin{enumerate}

\item{Every point $p \in \operatorname{States}$.}

\item{If $l \in \operatorname{States}$, then $l\overrightarrow{f} \in \operatorname{States}$.}

\end{enumerate}

\ed

In other words $\operatorname{States}$ is the smallest subset of $LV$ containing the points and closed under the action of $M_{0}(L)$ on $LV$. We can finally define the Evaluation Transformation Semigroup where all the computations in the examples in this paper take place.

\bd

The Evluation Transformation Semigroup $\mathcal{E}(L)$ is defined by $\mathcal{E}(L) = (\operatorname{States}, \operatorname{Eval}(L))$ where $\operatorname{Eval}(L)$ is the image of $M'_{0}(L)$ by restricting its action to $\operatorname{States}$.

\ed

We remark that there is an Evaluation Transformation Semigroup $\mathcal{E}_{n}(L)$ for all $n \geq 0$. The definition here is the case $n=0$
which is all we need in this paper as we are only concerned with aperiodic flows.

\subsection{Notation for Elements of the Rhodes Lattice $Rh_{B}(G)$}\label{SPCNotation}

In the examples  we compute using the Evaluation Transformation Semigroup $\mathcal{E}(L) = (\operatorname{States}, \operatorname{Eval}(L))$ where $L$ is the Rhodes lattice $Rh_{B}(G)$. We fix the notation we use for $SPC$ in our examples. Let $(W,\pi,[\mu]_\pi) \in Rh_B(G)$ be an element in the Rhodes lattice. Let $\pi_1,\pi_2,\ldots,\pi_k$ be the equivalence classes of
$\pi$ and let $\mu_i = \mu|_{\pi_i}: \pi_i \to G$ for $i\in \{1,2,\ldots,k\}$. Write $\pi_{1}\mid \pi_{2}\mid \ldots \pi_{k}/\langle \mu_1 \mid \mu_2 \mid \ldots \mid \mu_k\rangle$ for $(W,\pi,[\mu]_\pi)$. For readability, we drop set-brackets around elements of a partition class.
\be
\label{example.downstairs}
Let $G=\{1,-1\}=Z_2$, $W=\{1,2,3,4,5,6\}$, and
$\pi=\{1\mid  2 \ 5 \mid  3 \ 4 \ 6\}$. Let $\mu: \{1,2,3,4,5,6\} \to Z_2$ be defined by
$\mu(1)=\mu(3)=\mu(5)=1$ and $\mu(2)=\mu(4)=\mu(6)=-1$. Then we write $\pi=\{1 \mid  2 \ 5 \mid 3 \ 4 \ 6\}/\langle1 \mid \ -1 \ 1 \mid \ 1 \ -1 \ -1\rangle$.
\ee

%
%
$b \in Dom(f_{a})$ we have $bf_{a}=C(b,a)$. 

\bibliography{stubib}
\bibliographystyle{abbrv}


\end{document}